\documentclass[11pt]{article}
\usepackage{fancyhdr}
\usepackage[centertags]{amsmath}
\usepackage{amsfonts}
\usepackage{graphics}
\usepackage{graphicx}
\usepackage{xcolor}
\usepackage{amssymb}
\usepackage{amsthm}
\usepackage{newlfont}
\usepackage{url}
\usepackage{epsfig}
\usepackage{ifthen}
\usepackage{ifpdf}
\usepackage{rotating}
\usepackage{caption}
\usepackage{subcaption}
\usepackage{amsmath}
\usepackage[noadjust]{cite}

\input amssym.def
\input amssym
\input cyracc.def

\def\rn{\mathbb{R}^n}

\def\r{\mathbb{R}}

\newtheorem{lemma}{Lemma}[section]

\newtheorem{te}{Theorem}[section]

\newtheorem{ass}{Assumption}

\newcommand{\be}{\begin{equation}}
\newcommand{\ee}{\end{equation}}
\newcommand{\ba}{\begin{array}}
\newcommand{\ea}{\end{array}}
\newcommand{\bee}{\begin{eqnarray*}}
\newcommand{\eee}{\end{eqnarray*}}
\newcommand{\bea}{\begin{eqnarray}}
\newcommand{\eea}{\end{eqnarray}}

\newcommand{\N}{\mathcal{N}}

\newcommand{\F}{\mathcal{F}}

\begin{document}

\title{SLiSeS: Subsampled Line Search Spectral Gradient Method for  Finite Sums}

\author{Stefania Bellavia  \thanks{Department of Industrial Engineering, University of Florence,
		Viale Morgagni, 40/44, 50134 Florence, Italy, 
   e-mail: {\tt stefania.bellavia@unifi.it}.} \thanks{ INDAM-GNCS Research group, Roma, Italy.}
\and {Nata\v{s}a Kreji\'c  \thanks{Department of Mathematics and Informatics, Faculty of Sciences,  University of Novi Sad, Trg Dositeja Obradovi\'ca 4, 21000 Novi Sad, Serbia, e-mail: {\tt
natasak@uns.ac.rs}. }}
\and { Nata\v sa Krklec Jerinki\'c  \thanks{Department of Mathematics and Informatics, Faculty of Sciences,  University of Novi Sad, Trg Dositeja Obradovi\'ca 4, 21000 Novi Sad, Serbia, e-mail: {\tt
natasa.krklec@dmi.uns.ac.rs}. }}
 \and {Marcos Raydan \thanks{Center for Mathematics and Applications (NovaMath), FCT NOVA, 2829-516, Caparica,
 		Portugal, e-mail: {\tt m.raydan@fct.unl.pt}.  }}}

\maketitle

\begin{center}
{\it Dedicated to the memory of Andreas Griewank}
\end{center}

\begin{abstract} {The spectral gradient method is known to be a powerful low-cost tool for solving large-scale
	 optimization problems. In this paper, our goal is to exploit its advantages in the stochastic optimization
	 framework, especially in the case of mini-batch subsampling that is often used in big data settings.  To allow
	the spectral coefficient to properly explore the underlying approximate Hessian spectrum, we keep the same subsample for
	 {   a prefixed number of } iterations before subsampling again. We analyze the required algorithmic features and the conditions
	 for almost sure convergence, and present initial numerical results that show the advantages of the proposed
	  method. }

\end{abstract}

\section{Introduction}

We are concerned with the minimization of the average of finitely many  possibly nonconvex
smooth functions
\begin{equation} \label{prb}
\min_{x \in \rn}  f(x):= \frac{1}{N} \sum_{i=1}^{N} f_i(x),
\end{equation}
where each $f_i: \rn \rightarrow \r$ for $1\leq i \leq N$ is bounded from below and   $f_{low}$ denotes  the objective function's  lower bound.

Problem (\ref{prb}) has its origin in large-scale data analysis applications. For instance,
a wide variety of problems arising in machine learning involve computing an approximate minimizer
of the sum of a loss function over a large number of training examples, where there is a large amount
of redundancy between examples. In those cases, it is almost mandatory to employ stochastic iterative
methods that update the prediction model based on a relatively small randomly chosen subset (or sample)
 of the training data  \cite{bbkm,Bottou_et_al}.

 Let us denote  the full sample by  $\N=\{1,2,...,N\}$ and,  independently of the applied iterative method,
 let us denote the randomly chosen subsample at iteration $k$  by $\N_k \subseteq \N$ where $|\N_k |=S\ll N.$
The function estimator obtained by averaging the functions $f_i$ in $\N_k$ is given by
 $$f_{\N_k}(x)=\frac{1}{S} \sum_{i \in \N_k} f_i(x), $$
 and the associated  gradient  estimator is 
$$ \nabla f_{\N_k}(x)=\frac{1}{S} \sum_{i \in \N_k} \nabla f_i(x). $$

For solving  (\ref{prb}), stochastic gradient (SG) algorithms of the form
$$ x_{k+1} = x_k - \gamma_k \nabla f_{\N_k}(x_k), $$
exhibit a convenient low computational cost, and as a consequence have become popular and successful for solving large-scale
machine learning problems, for which the  full sample approach is prohibitively expensive. However, when combined with
standard choices of the step length $\gamma_k >0$ (constant or decreasing), these methods must perform a large number of iterations
to observe an adequate reduction in the objective function \cite{bkms,Bottou_et_al,cs,tmdq}.

Recently, for solving  (\ref{prb}), the SG method has been enriched with the use
 of spectral step lengths to speed up the convergence of the iterative process; see, e.g., \cite{bkm, kkj, sdr, tmdq, ywzl}.
 The spectral method (originally introduced by Barzilai and  Borwein \cite{BB88})
  and its variants, for the general (full sample) unconstrained minimization problem,  are low-cost gradient methods which have
proved to be very effective in practice for large-scale optimization. They have received a lot of attention in the last three decades
 including theoretical understanding, extensions, and adaptations for different scenarios (unconstrained and constrained)
  and for some specific applications; see, e.g., the review papers \cite{bmr14, Flet05, zoumag} and references therein. 
  The recently proposed schemes in  \cite{bkm, kkj, sdr, tmdq, ywzl}, which combine the SG approach with spectral step lengths, have a common feature: they all maintain
   the traditional option of changing the random subsample $\N_k$ (with constant or dynamically increasing size) at each iteration.

  From the very beginning it has been recognized, first for quadratics and later for the general unconstrained
   minimization problem, that the effectiveness of the spectral methods relies on the relationship between the inverse
   of the step lengths and the eigenvalues of the underlying Hessian of the objective function, near the local minimizer
    that is being approximated; see \cite{DaiFl, dhsz, ddrt13, drtz18, fhk23, Flet90, fzz08, ghr93, ray93}. The key issue that has been extensively 
    observed in practice is that, in sharp contrast
    with the classical Cauchy gradient method, the inverse of the step lengths of the spectral methods move all over the spectrum of
   the mentioned Hessian. As a consequence, these iterative methods exhibit a (highly) non-monotone behavior 
   during the convergence process.
   That unconventional behavior of the spectral gradient method, and its extensions, has been formally studied and 
   accurately called the {\it sweeping spectrum behavior}; see \cite{ddrt13, drtz18}.

  Now, back to solving (\ref{prb}), we note that changing the sample at each iteration, i.e., changing the chosen functions and the  gradients,
   produces a poor approximation of the inverse of the spectral step lengths to the  eigenvalues of the underlying Hessian
    due to the influence of noise. Therefore, changing the sample at each iteration is counterproductive to the effectiveness of
    the sweeping spectrum behavior.

   In order to promote the speed-up that can be obtained from the sweeping process,
   our main contribution in this work is to keep the same sample for {   a prefixed number of } iterations before subsampling again. That way, the inverse 
   of the step lengths can be focused for a few iterations on a set of fixed eigenvalues and enable the spectral method to show its 
   power in reducing the value of the objective function. 
    As we will show in our numerical experiments, this  change in the algorithm results in a significant improvement
     in the practical behavior of the method. Nevertheless, this novel strategy brings in new challenges with respect to the convergence analysis.

    The rest of this document is organized as follows.  In Section 2, we describe in detail the new algorithm 
     combining the new subsampling strategy with an appropriate line search globalization technique aimed at enhancing the practical behaviour of the overall procedure. In Section 3, we establish the convergence properties of the algorithm under some  variance assumptions for both uniform and nonuniform sampling.  We prove the almost sure convergence assuming 
     some additional properties 
         on the sampling. Moreover, we suggest a modification of the  algorithm which  employs a combination of predefined and adaptive step sizes and   allows to prove the convergence result avoiding both aforementioned additional assumptions. 
     
      In Section 4, we report on the obtained numerical results and  give further insights into the
    proposed approach. Finally, in Section 5 we present some concluding remarks.

\section{The SLiSeS Framework}

At each iteration, our scheme computes a mini-batch stochastic gradient estimate $g_k$ and employs diminishing step sizes based on the Barzilai-Borwein spectral coefficient.   
A nonmonotone line search is used that  employs  
 stochastic function estimates, at the current iterate and at the tentative new iterate, obtained by averaging the objective function in the same sample set used for the gradient estimator. 

Our  framework algorithm is described in Algorithm 1, while  Algorithm 2 details  the line search procedure. We will refer to the proposed framework as  SLiSeS (Subsampled Line Search Spectral gradient).

{\bf{Algorithm 1:  SLiSeS (Subsampled Line Search Spectral Gradient Method)}}

\begin{itemize}
\item[S0] Initialization: $x_0 \in \rn,  \eta \in (0,1), S \in \{1,2,...,N\},   m \in \mathbb{N}, 0< \gamma_{min}\leq 1 \leq \gamma_{max} < \infty, \{t_k\} \in \mathbb{R}^{\infty}_{+}$ such that $\sum_{k} t_k \leq \bar{t} < \infty$ and $ \N_0 \subseteq \N$, $|\N_0|=S$.
Set $k=0$.
\item[S1]  Sampling: If $mod(k,m)=0$ choose $ \N_k \subseteq \N$  such that $|\N_k |=S.$ Else, set $\N_k=\N_{k-1}$.
\item[S2]  Compute $g_k=\nabla f_{\N_k}(x_k)$.
\item[S3]  Spectral coefficient: If  $mod(k,m)=0$ and $m>1$ set $c_k=1/\|g_k\|$. Else, $ s_{k-1} = x_k - x_{k-1}, \; y_{k-1} = g_k - g_{k-1} $ and set $c_k=\|s_{k-1}\|^2 / (s_{k-1}^T y_{k-1})$.
\item[S4] Set $\gamma_k=\min \{\gamma_{max}, \max \{\gamma_{min},c_k\} \}/k$.
\item[S5] Search direction: Set $d_k=-\gamma_k \nabla f_{\N_k}(x_k).$
\item[S6] Step size: Find $\alpha_k \in (0,1]$ such that
$$f_{\N_k}(x_k+\alpha_k d_k)\leq f_{\N_k}(x_k)+\eta \alpha_k d^{T}_k \nabla f_{\N_k}(x_k)+t_k$$
by employing Algorithm 2.
\item[S7] Set $ x_{k+1}=x_k+\alpha_k d_k, \; k=k+1$ and go to Step S1.
\end{itemize}

 {At step S1 we compute the subsample used to form the gradient estimator at Step S2. This subsample is of size $S$ which is assumed to be significantly smaller than $N$ so that the mini-batch approach is employed at each iteration. 
Step S1 indicates that we keep the same subsample for $m$ iterations.  After every $m $ iterations, a new subsample is chosen. 
One of the main ingredients of the proposed algorithm is the spectral coefficient. At Step S3 of Algorithm 1 we calculate the coefficient $c_k$. It is  based on the spectral coefficient formula at the iterations where the subsample is not changed. On the contrary, after every $m$ iterations - when a new subsample $\N_k$ is introduced, $c_k$ is set to  $1/\|g_k\|$, except for the case  $m=1$. In this latter case the subsample is changed at each iteration and we 
 use the Spectral BB step length.  This way the method reduces to a stochastic gradient method with BB
 choice of the steplength similar to one considered in \cite{bkm}. 
 
 Furthermore, scaling with the gradient norm is relevant from practical point of view, while it is not crucial in the analysis and some other choices are feasible too, such as $c_k=1$ for instance. {   Occasional avoidance of the spectral coefficient is due to the fact that we let the subsample size $S$ to be arbitrary small and thus the noise can be rather large. In that case  the spectral coefficient might not offer any useful information on the spectrum of the Hessian as it would be computed on different estimator functions.}  Hence, we let the $c_k$ go through the spectrum of the approximate Hessian  only when the subsample is not changed, hoping that it will represent well its full sample counterpart. However, this may still provide some poor approximations given that we have a mini-batch instead of growing batch procedure. Moreover, we also aim to solve  nonconvex problems which can yield negative values of the spectral coefficients as well. Therefore, in step S4, we ``tame" the coefficient $c_k$ by projecting it on some positive, arbitrary large interval $[\gamma_{min}, \gamma_{max}]$ and by dividing it by $k$. Numerical experiments show that this strategy is beneficial since otherwise the algorithm exhibits erratic behavior.  After determining the  coefficient $\gamma_k$, we scale the negative gradient to obtain the search direction and use the line search procedure as described in Algorithm 2. 

Some comments on the Algorithm are in order.
\begin{itemize}
    \item Step S4 still allows enough freedom for the spectral coefficient; in the numerical result section we will show that SLiSeS outperforms stochastic gradient methods that employs  predefined sequence of step sizes such as $1/k$. 
\item The aforementioned noise of the spectral coefficient can be avoided by additional sampling. For instance, one can choose a new sample at each iteration and then calculate $\nabla f_{\N_k} (x_{k-1})$ to obtain $y_{k-1}=\nabla f_{\N_k} (x_{k})-\nabla f_{\N_k} (x_{k-1})$  and the corresponding $c_k$. Although this strategy is appealing, we believe that the approach presented in Algorithm 1 is more beneficial. The first reason is the computational cost which can be significantly enlarged if additional gradient sampling is applied, especially in high-dimensional problems. Second reason comes from the observed behavior of spectral method to collect the second order information by ``walking" through the spectrum which we allow by keeping the same sample for {   a prefixed number of } iterations. The resulting algorithm is rather cheap, efficient and supported by the underlying theoretical foundations presented in Section 3. 
\item We provide convergence analysis for both uniform and nonuniform sampling strategies, thus allowing biased estimators of the functions and the gradients as well. However, we assume that the approximate gradient is evaluated over the same subsample as the approximate function (Step S2). This allows us to obtain some valuable results concerning the line search procedure as indicated in Lemma \ref{lema1}. Although we allow nonmonotone line search in step S6 and thus any search direction is acceptable, S2 is crucial for proving that the step size is bounded away from zero under the assumption of Lipschitz-continuous gradients. This feature is sufficient for the convergence analysis that we perform in the sense that the  proofs are not leaning on the sufficient decrease condition coming from step S6. From theoretical point of view, $\alpha_k$ in Step S6 can even be a constant. 
However, as we will show in the experiments, the line search that we perform is highly relevant in practical behavior. 
\end{itemize}

 The Line Search procedure used to compute 
$\alpha_k$ in Step S6 of the SLiSeS Algorithm is described in 
Algorithm 2 and it is denoted as LSP (Line Search procedure).  At Step S3 of SLP the  quadratic interpolation provides a candidate step size $\tilde{\alpha}_j$;  if it  becomes too small, in Step S4 LSP switches to the backtracking procedure. The steplenght computed by LSP  remains bounded away from zero under some standard assumptions as it will be  shown in the next section. Although some other approaches are feasible, we choose this strategy motivated by our numerical experience. Notice that it does not require further sampling except for the evaluation of the function estimator  at new trial points needed to check the Armijo-like condition. 
\vskip 5 pt

\noindent {\bf{Algorithm 2:  LSP (Line Search Procedure)}}
\begin{itemize}
\item[S0] Input parameters: $x_k \in \rn,  \eta \in (0,1), \N_k, \gamma_k, t_k$.
\item[S1] Initialization: Set $j=0, \; dm_k=g_k^T d_k, \;  \alpha_j=1$.
\item[S2] If
\begin{equation} \label{Armijo}
    f_{\N_k}(x_k+\alpha_j d_k)\leq f_{\N_k}(x_k)+\eta \alpha_j dm_k+t_k
\end{equation}
go to Step S5. Else, go to Step S3.
\item[S3] If $\alpha_j >  0.1,$ compute $$\tilde{\alpha}_j=\frac{-dm_k \alpha_j^2}{2 (f_{\N_k}(x_k+\alpha_j d_k)-f_{\N_k}(x_k)-\alpha_j dm_k)}.$$
 If $\tilde{\alpha}_j < 0.1 \alpha_j$ or $\tilde{\alpha}_j>0.9 \alpha_j$, set $\tilde{\alpha}_j = \alpha_j /2$. Set $\alpha_{j+1}=\tilde{\alpha}_j, \; j=j+1$ and  go to Step S2.
 \item[S4]  If $\alpha_j\leq 0.1$, set $\alpha_{j+1}=\alpha_j /2, \; j=j+1$ and  go to Step S2.
\item[S5] Set $\alpha_k=\alpha_j$ and STOP.
\end{itemize}

\section{Convergence analysis}

{ Within this section we analyze conditions under which the spectral gradient method SLiSeS converges a.s. in stochastic mini-batch framework. We start our analysis by proving that LSP yields uniformly bounded step sizes under the assumption of Lipschitz-continuous gradients stated in A\ref{A1}.
Then, Subsection 3.1 considers uniform sampling which yields unbiased estimators. 
Assuming the common  bounded variance assumption A\ref{A2}, we prove a.s. convergence result under 
the Assumption A\ref{A2ss}  on the sampling properties (Subsection 3.1.1). We also prove that the a.s. convergence result can be achieved by avoiding Assumption  A\ref{A2ss}, provided that  the proposed algorithm is modified as outlined  in Subsection {   3.1.2.} 
{   Subsection 3.2 is devoted to nonuniform sampling as stated in Assumption A\ref{A5}. Under the similarity assumption stated in A\ref{A4}, we prove the a.s. convergence for 
the modified  algorithm as in Subsection    3.1.2.} Finally, we show that stronger convergence result can be achieved for {   functions satisfying the Polyak- Lojasiewicz (P-L) condition} and for strongly convex functions in Subsection 3.3.  }

\begin{ass}[\bf{A1}] \label{A1} 
 The functions $f_i, i=1,...,N$ are bounded from below and  twice   continuously-differentiable   with  L-Lipschitz-continuous gradients.
\end{ass}

Notice that Assumption A\ref{A1} implies that $ \|\nabla^2 f_i(x) \|\leq L$. Moreover, for any $\N_k$ the function $f_{\N_k}$ is bounded from below with L-Lipschitz-continuous gradients and bounded Hessians, i.e., $ \|\nabla^2 f_{\N_k}(x)\|\leq L.$

As already mentioned, the line search procedure is well defined. This is shown in the following Lemma.
\begin{lemma} \label{lema1} Suppose that A\ref{A1} holds. Then the LSP procedure is well defined and there exists a constant  $\alpha_{min}>0 $ such that $\alpha_{min} \leq \alpha_k\leq 1 $ for every $k$.
\end{lemma}

 \begin{prf} Since the function $f_{\N_k}$ is bounded from below, and the direction $d_k$ is descent,  there exists an interval of step sizes which satisfy Armijo condition  \eqref{Armijo}.

  Let us now show that the step sizes are bounded from below. The step size is computed either in S3 or S4 of LSP. Assume first that the current value of $ \alpha_j  $ is such that $ \alpha_j > 0.1$ so we enter S3 to compute $\tilde{\alpha}_j. $ If $ \tilde{\alpha}_j <0.1 \alpha_j $ then we have $ \alpha_{j+1} = \alpha_j/2 > 0.05. $ In all other cases we have $ \alpha_{j+1} > 0.1 \alpha_j > 0.01.  $ Thus for all $ \alpha_j $ computed in S3 we have $ \alpha_j > 0.01. $ Thus we can enter S4 with $ \alpha_j < 0.1 $ but at the same time $ \alpha_j > 0.01,$ i.e. the value that is first used in S4 can not be arbitrary close to 0. Observe that once we enter S4 we can not go back to S3, and  all subsequent values of the step size are computed by S4. Let us denote by
 $\alpha'_k$ the last step size in LSP which did not satisfy \eqref{Armijo}, i.e., $\alpha'_k=2 \alpha_k$. Since $t_k\geq 0$ there holds
 $$f_{\N_k}(x_k+\alpha'_k d_k)>  f_{\N_k}(x_k)+\eta \alpha'_k d^{T}_k \nabla f_{\N_k}(x_k).$$
 On the other hand, using the Lipschitz continuity  we obtain
 $$ f_{\N_k}(x_{k}+\alpha'_k d_k)    \leq f_{\N_k}(x_k)  + \alpha'_k \nabla f_{\N_k}(x_k)^T d_k + \frac{L}{2} {\alpha'}_k^2 \|d_k\|^2.  $$
 Combining the previous two inequalities and
$ d_k = -\gamma_k \nabla f_{\N_k}(x_k), $ with $ \gamma_k \leq \gamma_{\max} $ as specified Step S4 of SLiSeS algorithm,   we obtain
 $$\alpha'_k \geq \frac{-2(1-\eta) \nabla f_{\N_k}(x_k)^T d_k}{L\|d_k\|^2}=\frac{2(1-\eta) \gamma_k \|\nabla f_{\N_k}(x_k)\|^2 }{L\gamma^2_k \|\nabla f_{\N_k}(x_k)\|^2}\geq \frac{2(1-\eta)}{L\gamma_{max}}.$$
 Thus,   $\alpha_k \geq (1-\eta)/(L\gamma_{max})$ for $ \alpha_k $ computed in S4. Recall that we already concluded that $ \alpha_k > 0.01 $ if it is computed in S3. Therefore the statement holds with
$$\alpha_{min}=\min \{\frac{(1-\eta)}{L\gamma_{max}}, 0.01\}.$$ $\Box$
 \end{prf}

Denote by $\F_k$ a $\sigma$-algebra generated by $\N_0, ...,\N_{k-1}$. Notice that $x_k$ is $\F_k$-measurable, i.e., $x_k$ is fully determined (known) under condition of knowing $\F_k$. On the other hand, $\alpha_k, g_k$ and $\gamma_k$ are not $\F_k$-measurable in general.  We will be interested only in $\F_k$ for iterations $k$ such that $mod(k,m)=0$. We will refer to these iterations as outer iterations. In the remaining iterations, to be referred to as inner iterations, we keep the same sample.

The following theorem from \cite{robbins:1971} is used for proving the a.s. convergence result.
\begin{te} \label{thm:rand_variable_conv}
    Let $ U_l, \beta_l, \xi_l, \rho_l \geq 0 $ be $ \F_l$-measurable random variables such that
    $$ E[U_{l+1} | \F_l] \leq  (1+\beta_l) U_l + \xi_l - \rho_l, \quad l=1,2,\ldots \, . $$
    If $ \displaystyle \sum_l \beta_l < \infty $ and $ \displaystyle \sum_l \xi_l < \infty $, then $ U_l \to U < \infty $ a.s. and $ \displaystyle \sum_l \rho_l < \infty $ a.s.
\end{te}

\subsection{Uniform sampling}

Within this subsection we assume that the sampling in step S1 of SLiSeS method is uniform. This implies that
\begin{equation} \label{unb}
E[f_{\N_k}(x_k) | \F_k]=f(x_k) \quad \mbox{and}  \quad E[ g_k | \F_k]= \nabla f(x_k),
\end{equation}
{   for all $ k $ such that $ mod(k,m)=0. $ } 

We state the following common assumption for subsampling methods.

\begin{ass}[\bf{A2}] \label{A2}  There exists a constant $G>0$ such that for all $k$ {   with  $ mod(k,m)=0 $}  there holds $$E[ \|g_k-\nabla f(x_k)\|^2 | \F_k]\leq G.$$
\end{ass}
Assumption A\ref{A2} together with \eqref{unb} implies
\begin{equation}
\label{gkG}
E[ \|g_k \|^2 | \F_k]\leq G+ \|\nabla f(x_k)\|^2,
\end{equation}
{   at any iteration {   $ k $ such that $ mod(k,m)=0. $}}

\subsubsection{Subsampling strategies }
 Given that $ \N_k$ is finite let $ s $ denote the   number of possible samples $ \N_k: \N_k^1,\ldots,\N_k^s. $ Then, for any $ \N_k $ {   at iteration  $ k $ such that $ mod(k,m)=0$}, we  have $\|\nabla f(x_k)\|^2 = E[\nabla f(x_k)^T\nabla f_{\N_k}(x_k) |\F_k] $ and
\begin{eqnarray} \label{novo1}
\|\nabla f(x_k)\|^2 = \frac{1}{s} \left(\sum_{i \in I_k} \nabla f(x_k)^T \nabla f_{\N_k^i}(x_k) + \sum_{j \in J_k} \nabla f(x_k)^T \nabla f_{\N_k^j}(x_k) \right),
\end{eqnarray}
where $ I_k =\{i:\nabla f(x_k)^T \nabla f_{\N_k^i}(x_k)\geq 0\} $ and  $ J_k =\{j:\nabla f(x_k)^T \nabla f_{\N_k^j}(x_k)< 0\}. $
 Notice that (\ref{novo1}) implies
that $ I_k \neq \emptyset $ for all $k$.  We state the following assumption.

 \begin{ass}[\bf{A3}] \label{A2ss}  There exist constants $C, \theta>0$ and $ 0< \delta <1 $ such that
 {   at any iteration  $ k $ such that $ mod(k,m)=0,$}
\begin{equation} \label{A2ss_ineq}
P(B_k) E[ \nabla f(x_k)^T \nabla f_{\N_k}(x_k)| \F_k, B_k] \leq \frac{\theta}{s}\|\nabla f(x_k)\|^2 + \frac{C}{s k^\delta},
\end{equation}
where $B_k $ is the event that  $\nabla f(x_k)^T \nabla f_{\N_k}(x_k) \geq 0$ and $P(B_k)$ is the probability of this event.
\end{ass}
  Notice that \eqref{A2ss_ineq} can also be obtained by assuming the following 
\begin{equation} \label{ass_trish}
E[||\nabla f_{\N_k}(x_k)||^2|  \F_k]\le \frac{\tilde{C}}{k^{\tilde{\delta}}}+\tilde{\theta }\|\nabla f(x_k)\|^2, 
\end{equation}
for some $\tilde{\theta}, \tilde{C}, \tilde{\delta} >0$ (see \cite[Lemma 2]{trish}). Thus, the Assumption A\ref{A2ss} can be replaced with \eqref{ass_trish}, but we decide to keep the weaker assumption that states \eqref{A2ss_ineq}. {   The assumption is closely related to the so-called orthogonality test in \cite{nocedal}, which is used as the sample increase test. Assuming that A3 holds could be interpreted as an assumption that the subsample is large enough, i.e. in our case that $ S $ is large enough.  Another example of functions that satisfy \eqref{ass_trish} with $\tilde C=0$
 are the ones that satisfy the following strong growth condition (SGC):
$$
\frac{1}{N} \sum_{i=1}^N||\nabla f_i(x_k)||^2 \le \tilde{\theta }\|\nabla f(x_k)\|^2,
$$
for a given $\tilde{\theta }$, coupled with uniform sampling.
Such functions arise in  
binary classification with a linear classifier employing the squared hinge loss or 
logistic regression  when the data is linearly separable, as shown in \cite[Appendix A]{Chen_et_al}.
}

\begin{te} \label{Tuniss}
  Suppose that the Assumptions A\ref{A1},A\ref{A2} and A\ref{A2ss} hold. Then, provided that $s\geq \theta$, with probability 1  there holds
$\liminf_{k \to \infty}  \|\nabla f(x_k)\|=0.$
\end{te}

\begin{prf}
{   Let us consider a subsequence $ K \subset \mathbb{N} $ such that $ k \in K $ if $ mod(k,m)=0. $ For the sake of readability we derive the proof for $ m=3$ but the reasoning is exactly the same for larger values of $m. $
So, given $ k \in K $ and $ x_k $,  we choose a sample $ \N_k $ and compute $ x_{k+1}$, $x_{k+2} $ and $ x_{k+3}$ taking the directions $ -\nabla f_{\N_k}(x_k)$,  $-\nabla f_{\N_k}(x_{k+1})$ and $ -\nabla f_{\N_k}(x_{k+2}) $ and applying the line search algorithm. Therefore we have
 $$ x_{k+3} = x_k + \sum_{j=0}^{2} \alpha_{k+j}  d_{k+j}$$ with
$ d_{k+j} = -\gamma_{k+j} \nabla f_{\N_k}(x_{k+j})$.
 Letting $\nabla^2 f_{\N_k}(x)$ be the Hessian matrix of $f_{\N_k}$ at $x$, 
the Taylor expansion yields
\begin{equation} \label{eq0}
  \nabla f_{\N_k}(x_{k+1}) =   \nabla f_{\N_k}(x_{k}) - \alpha_k \gamma_k \nabla^2 f_{\N_k}(\theta_k) \nabla f_{\N_k}(x_{k}),
\end{equation}
for some $ \theta_k \in B(x_k,\alpha_k\|d_k\|),$ and
\begin{eqnarray} \label{eq1}
  \nabla f_{\N_k}(x_{k+2}) &  = &  \nabla f_{\N_k}(x_{k+1}) - \alpha_{k+1}\gamma_{k+1} \nabla^2 f_{\N_k}(\theta_{k+1}) \nabla f_{\N_k}(x_{k+1}) \nonumber \\
  & = &   \nabla f_{\N_k}(x_{k}) - \alpha_k \gamma_k \nabla^2 f_{\N_k}(\theta_k) \nabla f_{\N_k}(x_{k}) -\nonumber \\
  && \alpha_{k+1}\gamma_{k+1} \nabla^2 f_{\N_k}(\theta_{k+1})(  \nabla f_{\N_k}(x_{k}) - \alpha_k \gamma_k \nabla^2 f_{\N_k}(\theta_k) \nabla f_{\N_k}(x_{k})) \nonumber \\
  & = &  \nabla f_{\N_k}(x_{k}) - \alpha_k \gamma_k \nabla^2 f_{\N_k}(\theta_k) \nabla f_{\N_k}(x_{k}) - \alpha_{k+1} \gamma_{k+1} \nabla^2 f_{\N_k}(\theta_{k+1}) \nabla f_{\N_k}(x_k)+ \nonumber \\
  && \alpha_k \gamma_k \alpha_{k+1}\gamma_{k+1} \nabla^2 f_{\N_k}(\theta_{k+1}) \nabla^2 f_{\N_k}(\theta_{k}) \nabla f_{\N_k}(x_k),
  \end{eqnarray}
  for some $ \theta_{k+1} \in B(x_{k+1}, \alpha_{k+1} \|d_{k+1}\|).$  
   Since $\alpha_k\leq1$ and $\gamma_k\leq \gamma_{max}/k\leq \gamma_{max}$ for all $k$,  the assumption A\ref{A1} further implies the existence of constants $C_1, C_2$ such that
\begin{eqnarray}
    \label{gk12}
    \|\nabla f_{\N_k}(x_{k+1})\|\leq C_1 \|\nabla f_{\N_k}(x_{k})\|, \quad \|\nabla f_{\N_k}(x_{k+2})\|\leq C_2 \|\nabla f_{\N_k}(x_{k})\|.
\end{eqnarray}
Furthermore, }
 we can express $ \nabla f(x_k)^T(x_{k+3}-x_k)$ as follows,
 \begin{eqnarray} \label{eq1a}
  \nabla f(x_k)^T(x_{k+3}-x_k)&=& - (\sum_{j=0}^2 \alpha_{k+j}\gamma_{k+j}) \nabla f(x_k)^T \nabla f_{\N_k}(x_k) +
 \nonumber \\
 &&  \alpha_{k+1}\gamma_{k+1} \alpha_{k}\gamma_{k}  \nabla f(x_k)^T \nabla^2 f_{\N_k}(\theta_k)\nabla f_{\N_k}(x_k) + \nonumber \\
 && \alpha_k \gamma_k \alpha_{k+2} \gamma_{k+2} \nabla f(x_k)^T \nabla^2 f_{\N_k}(\theta_k)\nabla f_{\N_k}(x_k) + \nonumber \\
&& \alpha_{k+1} \gamma_{k+1} \alpha_{k+2} \gamma_{k+2} \nabla f(x_k)^T \nabla^2 f_{\N_k}(\theta_{k+1})\nabla f_{\N_k}(x_k)- \nonumber \\
&&\alpha_k \gamma_k \alpha_{k+1} \gamma_{k+1} \alpha_{k+2} \gamma_{k+2} \nabla f(x_k)^T \nabla^2 f_{\N_k}(\theta_{k+1})\nabla^2 f_{\N_k}(\theta_k) \nabla f_{\N_k}(x_k).
 \end{eqnarray}
{    Notice}
  that for all   $ k \in \mathbb{N}$ we have
 \begin{equation}\label{eqstepsize}
  \frac{ \alpha_{\min} \gamma_{\min}}{k} \leq \alpha_k \gamma_k \leq \frac{\gamma_{\max} }{k}.
 \end{equation}
 Moreover, since the Assumption A\ref{A1} implies bounded Hessians, we obtain
\begin{eqnarray} \label{eq8}
   \nabla f(x_k)^T(x_{k+3}-x_k)  &\leq&   - (\sum_{j=0}^2 \alpha_{k+j}\gamma_{k+j}) \nabla f(x_k)^T \nabla f_{\N_k}(x_k) +\nonumber \\
  & & \frac{\gamma^2_{\max}}{k^2}L\|\nabla f(x_k)\| \|\nabla f_{\N_k}(x_k)\|(3 + \frac{\gamma_{\max}}{k}L) \nonumber  \\
  & =  & - (\sum_{j=0}^2 \alpha_{k+j}\gamma_{k+j}) \nabla f(x_k)^T \nabla f_{\N_k}(x_k) + \frac{\tilde{C}_1}{k^2} \|\nabla f(x_k)\| \|\nabla f_{\N_k}(x_k)\|,
\end{eqnarray}
with $ \tilde{C}_1 = \gamma_{\max}^2L(3 + \gamma_{\max}L).$
Furthermore, applying \eqref{gk12} we obtain 
\begin{eqnarray} \label{eq8a}
    \|x_{k+3}-x_k\|^2  &\leq&  4 \sum_{j=0}^2 \|\alpha_{k+j} d_{k+j}\|^2\le 4 \frac{\gamma_{\max}^2}{k^2} \|\nabla f_{\N_k}(x_k)\|^2 \big(1+C_1^2+C_2^2\big )
    \nonumber \\
    & = & \frac{\tilde{C}_2}{k^2}  \|\nabla f_{\N_k}(x_k)\|^2,
\end{eqnarray}
with $ \tilde{C}_2 = 4 \gamma_{\max}^2(1+C_1^2+C_2^2).$
Using Assumption A\ref{A1} and  the descent lemma we get
$$ f(x_{k+3})    \leq f(x_k)  +\nabla f(x_k)^T(x_{k+3}-x_k) + \frac{L}{2} \|x_{k+3}-x_k\|^2.  $$
 Subtracting $f_{low}$ and applying  the conditional expectation with respect to $\F_k$ we obtain
\begin{equation} \label{p1}
        E[f( x_{k+3}) -f_{low}| \F_k] \leq f( x_{k}) -f_{low}+ E[\nabla f(x_k)^T(x_{k+3}-x_k)| \F_k] + \frac{L}{2}  E[\|x_{k+3}-x_k\|^2 | \F_k].
\end{equation}
\noindent
{   Further, by \eqref{gkG}} and \eqref{eq8a} we  obtain
\begin{eqnarray} \label{p22}
E[\|x_{k+3}-x_k\|^2 | \F_k]\leq \frac{\tilde{C}_2}{k^2} E[\|\nabla f_{\N_k}(x_k)\|^2 | \F_k] \leq {   \frac{ \tilde{C}_2 }{k^2} (G+\| \nabla f(x_k)\|^2)}.
\end{eqnarray}
Now, let us upper bound  $E[\nabla f(x_k)^T(x_{k+3}-x_k) | \F_k] $.  Starting from \eqref{eq8} and using the fact that $\nabla f(x_k)$ is $\F_k$-measurable we get
\begin{eqnarray} \label{eq11}
     E[\nabla f(x_k)^T(x_{k+3}-x_k) | \F_k] & =& E[-\sum_{j=0}^2 \alpha_{k+j} \gamma_{k+j}\nabla f(x_k)^T \nabla f_{\N_k}(x_k) | \F_k] +\nonumber \\
    & & E[\frac{\tilde{C}_1}{k^2}\|\nabla f(x_k)\| \|\nabla f_{\N_k}(x_k)\| |\F_k] \nonumber\\
    & = & E[-\sum_{j=0}^2 \alpha_{k+j} \gamma_{k+j}\nabla f(x_k)^T \nabla f_{\N_k}(x_k) | \F_k] + \nonumber \\
    && \frac{\tilde{C}_1}{k^2}\|\nabla f(x_k)\| E[\|\nabla f_{\N_k}(x_k)\| |\F_k].
\end{eqnarray}
 Let us upper bound the second term first. There holds
\begin{equation}\label{cor5}E[\|\nabla f_{\N_k}(x_k)\| |\F_k]\leq \sqrt{E[\|\nabla f_{\N_k}(x_k)\|^2 |\F_k]}\\
\leq     {\sqrt{G} +\|\nabla f (x_k)\|}.
\end{equation}
 Thus, considering both cases, i.e., $\|\nabla f (x_k)\|\leq 1$ and $\|\nabla f (x_k)\|> 1$, we obtain
\begin{equation}\label{novo2}
\frac{\tilde{C}_1}{k^2}\|\nabla f(x_k)\| E[\|\nabla f_{\N_k}(x_k)\| |\F_k]\leq
 \frac{\tilde{C}_3}{k^2}{  (\sqrt{G}+\|\nabla f(x_k)\|^2)},
 \end{equation}
 where {   $\tilde{C}_3:=\tilde{C}_1(\sqrt{G}+1).$}
Finally, to upper bound $E[-\sum_{j=0}^2 \alpha_{k+j} \gamma_{k+j}\nabla f(x_k)^T \nabla f_{\N_k}(x_k) | \F_k] $ we proceed as follows.
For any realizations such that $B_k$ holds,
using (\ref{eqstepsize}) we have
\begin{equation} \label{Bktrue}
-\sum_{j=0}^2 \alpha_{k+j} \gamma_{k+j}\nabla f(x_k)^T \nabla f_{\N_k}(x_k)\le - 3 \frac{\gamma_{\min}  \alpha_{min}}{k+2} \nabla f(x_k)^T \nabla f_{\N_k}(x_k),
\end{equation}
while if $B_k$ is not true
\begin{equation} \label{Bkfalse}
-\sum_{j=0}^2 \alpha_{k+j} \gamma_{k+j}\nabla f(x_k)^T \nabla f_{\N_k}(x_k)\le - 3 \frac{\gamma_{\max}}{k} \nabla f(x_k)^T \nabla f_{\N_k}(x_k).
\end{equation}
Moreover, letting $\overline B_k$ the complementary event of $B_k$ it holds:
\begin{eqnarray}
E[-\sum_{j=0}^2 \alpha_{k+j} \gamma_{k+j}\nabla f(x_k)^T \nabla f_{\N_k}(x_k) | \F_k]&=&
P(B_k) E[-\sum_{j=0}^2 \alpha_{k+j} \gamma_{k+j}\nabla f(x_k)^T \nabla f_{\N_k}(x_k) | \F_k,B_k]+\nonumber\\
&& P(\overline B_k) E[-\sum_{j=0}^2 \alpha_{k+j} \gamma_{k+j}\nabla f(x_k)^T \nabla f_{\N_k}(x_k) | \F_k,\overline B_k] - \nonumber\\
&\le& 
 {\color{blue} -} 3 \frac{\gamma_{\min}  \alpha_{min}}{k+2} P(B_k)E[\nabla f(x_k)^T \nabla f_{\N_k}(x_k) | \F_k,B_k] \nonumber\\
 &-&  3 \frac{\gamma_{\max}}{k} P(\overline B_k)E[\nabla f(x_k)^T \nabla f_{\N_k}(x_k) | \F_k,\overline B_k]. \label{E1}
 \end{eqnarray}

We further notice that 
\begin{eqnarray*}
{  \|\nabla f(x_k)\|^2}&=&E[\nabla f(x_k)^T \nabla f_{\N_k}(x_k) | \F_k] \\
&=& P(B_k)E[\nabla f(x_k)^T \nabla f_{\N_k}(x_k) | \F_k,B_k]+P(\overline B_k)E[\nabla f(x_k)^T \nabla f_{\N_k}(x_k) | \F_k,\overline B_k].
\end{eqnarray*}
Then,
$$
-P(\overline B_k)E[\nabla f(x_k)^T \nabla f_{\N_k}(x_k) | \F_k,\overline B_k]=- {   \|\nabla f(x_k)\|^2}+ P(B_k)E[\nabla f(x_k)^T \nabla f_{\N_k}(x_k) | \F_k,B_k],
$$
and by \eqref{E1} and Assumption A\ref{A2ss}  we get
\begin{eqnarray}
E[-\sum_{j=0}^2 \alpha_{k+j} \gamma_{k+j}\nabla f(x_k)^T \nabla f_{\N_k}(x_k) | \F_k]  &{   \le}&
 - 3 \frac{\gamma_{\min}  \alpha_{min}}{k+2} P(B_k)E[\nabla f(x_k)^T \nabla f_{\N_k}(x_k) | \F_k,B_k]\nonumber\\
 &&+ 3 \frac{\gamma_{\max}}{k} 
 (- {   \|\nabla f(x_k)\|^2}+ P(B_k)E[\nabla f(x_k)^T \nabla f_{\N_k}(x_k) | \F_k,B_k])\nonumber \\
 &=&
 3 (P(B_k)E[\nabla f(x_k)^T \nabla f_{\N_k}(x_k) | \F_k,B_k] (\frac{\gamma_{\max}}{k} -\frac{\gamma_{\min}  \alpha_{min}}{k+2})\nonumber\\
 &&-
  \frac{\gamma_{\max}}{k} {   \|\nabla f(x_k)\|^2})\nonumber\\
 &\le & 
3 ((\frac{\gamma_{\max}}{k} -\frac{\gamma_{\min}  \alpha_{min}}{k+2}) 
( \frac{\theta}{s} \|\nabla f(x_k)\|^2 + \frac{C}{s  k^\delta})-
  \frac{\gamma_{\max}}{k} {   \|\nabla f(x_k)\|^2} )\nonumber\\
  &\le&
  3 (-{   \|\nabla f(x_k)\|^2}(\frac{\gamma_{\max}}{k}+\frac{\theta \gamma_{\min}  \alpha_{min}}{s(k+2)}-
 \frac{\theta \gamma_{\max}}{s k} )+ 
  \frac{ \gamma_{\max} C}{s  k^{1+\delta}})\nonumber\\
  &\le&
  - \frac{\tilde C_4}{k+2}{   \|\nabla f(x_k)\|^2}+ 
  \frac{\tilde C_5}{  k^{1+\delta}} \label{novo3}
\end{eqnarray}
where   ${  \tilde{C}_4:= 3 (\theta\gamma_{min} \alpha_{min} / s)}$, $\tilde{C}_5:=3 \gamma_{max} C/s$ and the last inequality holds as $s\ge \theta$ by assumption.

Finally, combining (\ref{p1}), (\ref{p22}), (\ref{novo2}) and (\ref{novo3}) we obtain
\begin{eqnarray} \label{s4}
        E[f( x_{k+3}) -f_{low} | \F_k) & \leq & f( x_{k}) -f_{low}-
         \left(\frac{\tilde{C}_4}{k+2} -\frac{\tilde{C}_3+{   \tilde{C}_2L/2}}{ k^2}\right) \|\nabla f(x_k)\|^2  +\frac{\tilde{C}_6}{k^{1+\delta}}\\\nonumber
        & = & f( x_{k}) -f_{low} -\frac{1}{k}\left( \frac{\tilde{C}_4 k}{k+2}-\frac{\tilde{C}_3+\tilde{C}_2  L/2}{k}\right)\|\nabla f(x_k)\|^2 +\frac{\tilde{C}_6}{k^{1+\delta}},
\end{eqnarray}
where $\tilde{C}_6:=\tilde{C}_3{   \sqrt{G}}+\tilde{C}_5+\tilde{C}_2  GL/2$.
 Since $\frac{\tilde{C}_4 k}{k+2} \to \tilde{C}_4$ and $\frac{\tilde{C}_3+\tilde{C}_2 L/2}{k}\to 0$ when $k\to \infty$, there exists $\bar{k} \in K$ such that for  all $ k \geq \bar{k}, k \in K $ we have
\begin{equation} \label{s5}
        E[f( x_{k+3}) -f_{low} | \F_k] \leq f( x_{k}) -f_{low} - \frac{{\bar{C}}_4}{2 k} \|\nabla f(x_k)\|^2 + \frac{{\bar{C}}_6}{k^{1+\delta}}.
\end{equation}
{   Thus, we have just proved that
  \begin{equation} \label{p6l} E[U_{\ell+1} | \F_\ell] \leq (1+\beta_\ell) U_l+\xi_\ell-\rho_\ell,
 \end{equation}
 holds for  $\ell=0,1,2,...$  where
 $$ U_{\ell+1}:= f( x_{\bar{k}+3(\ell+1)}) -f_{low}, \;
  \xi_\ell:=\frac{{\bar{C}}_6 }{(\bar{k}+3\ell)^{1+\delta}},\;  \rho_\ell:=\frac{\bar{C}_4}{2(\bar{k}+3\ell)} \|\nabla f(x_{\bar{k}+3\ell})\|^2,  \; \beta_l=0. $$
 Notice that $U_\ell, \xi_\ell, \rho_\ell, \beta_\ell$ are nonnegative and $\F_\ell$-measurable. Moreover,  there holds $\sum_\ell \xi_\ell < \infty$ and obviously $\sum_\ell \beta_\ell < \infty$. Thus, applying Theorem~\ref{thm:rand_variable_conv} we conclude that
    $$ \sum_{\ell=1}^{\infty} \rho_\ell= \frac{\bar{C}_4}{2 }\sum_{\ell=1}^{\infty} \frac{\|\nabla f(x_{\bar{k}+3\ell})\|^2}{\bar{k}+3\ell} < \infty \quad \mbox{a.s.} $$
    This further implies that $\liminf_{\ell \to \infty} \|\nabla f(x_{\bar{k}+3\ell})\|^2 =0$ a.s., which completes the proof.}
     $\Box$
\end{prf}

\subsubsection{Algorithmic modifications}
We have proved the a.s. convergence result under assumptions A\ref{A1}, A\ref{A2} and  an additional assumption 
on the sampling strategy (A\ref{A2ss}).  Focusing on the case $m>1$, 
we can avoid additional assumptions and retain the a.s  convergence by slightly modifying the steps S3-S5 of SLiSeS algorithm   as follows. {   Take $ \delta > 0,$ arbitrary small and $ \tilde{\gamma} > 0$. }
\begin{itemize}
\item[S3']  If  $mod(k,m)=0$  set  $\gamma_k=\tilde{\gamma}/k$ and go to step S5'. Else, set  $ s_{k-1} = x_k - x_{k-1}, \; y_{k-1} = g_k - g_{k-1} $, set $c_k=\|s_{k-1}\|^2 / (s_{k-1}^T y_{k-1})$  and go to step S4'.
\item[S4']  $\gamma_k=\min \{\gamma_{max}, \max \{\gamma_{min},c_k\} \}/k^{1+\delta}$.
\item[S5'] Set $d_k=-\gamma_k \nabla f_{\N_k}(x_k).$  If $mod(k,m)=0$, set $\alpha_k=1$ and go to step S7. Else, go to step S6.
\end{itemize}
 {   We note that the choice $\gamma_k=\tilde{\gamma}/k$ when $mod(k,m)=0$ yields $\F_k$-measurable steplengths whenever the subsample is recomputed {   as $\gamma_k$ and $\alpha_k$ are not dependent of the subsample ${\cal N}_k$}. Notice that taking $\gamma_k=\tilde{\gamma}/k$ instead of $\gamma_k=1/k$ in step S3 does not change the theoretical analysis and for simplicity we will consider $ \tilde{\gamma}=1 $ in the further analysis. The choice of $ \tilde{\gamma} $ might influence the numerical performance of the algorithm and in our experiments we used $\tilde{\gamma}=1/\|g_0\|$ where $g_0$ is  a stochastic estimator of the gradient at $x_0$, for further details see Section 4. } 
 We further underline that  we use  stepsizes diminishing as   $O(\frac{1}{k^{1+\delta}})$ rather than  as $O(\frac{1}{k})$. These two conditions allow to  prove convergence without relying  on 
 A\ref{A2ss} as shown in the following theorem. 
\begin{te} \label{Tunimodi} Suppose that the Assumptions A\ref{A1} and A\ref{A2} hold. Let $\{x_k\}$ be a sequence generated by  modified SLiSeS algorithm with  
$m>1$ and steps S3'-S5' in place of  of S3-S5. Then, with probability 1, there holds
$\liminf_{k \to \infty}  \|\nabla f(x_k)\|=0.$
\end{te}

\begin{prf}
{   For simplicity we take $ \tilde{\gamma}=1 $ in the proof. Following the same steps as in Theorem \ref{Tuniss} until \eqref{novo2}, we obtain \eqref{eq8},
\eqref{p1} and \eqref{p22}. 
Then, to upper bound $E[-\sum_{j=0}^2 \alpha_{k+j} \gamma_{k+j}\nabla f(x_k)^T \nabla f_{\N_k}(x_k) | \F_k] $, notice that $\gamma_k=1/k$ and  $\alpha_k=1$, so they   are {   not dependent on the subsample  $\cal{N}_k$ and hence they are }$\F_k$-measurable. Thus,
\begin{eqnarray} \label{sl3aaa} 
 E[-\sum_{j=0}^2 \alpha_{k+j} \gamma_{k+j}\nabla f(x_k)^T \nabla f_{\N_k}(x_k) | \F_k] & = & - \frac{1}{k} \nabla f(x_k)^T E[\nabla f_{\N_k}(x_k) | \F_k]\nonumber \\
&& -  E[\sum_{j=1}^2 \alpha_{k+j} \gamma_{k+j}\nabla f(x_k)^T \nabla f_{\N_k}(x_k) | \F_k]
\end{eqnarray}
Moreover, there holds 
$$-\sum_{j=1}^2 \alpha_{k+j} \gamma_{k+j}\nabla f(x_k)^T \nabla f_{\N_k}(x_k) \leq \sum_{j=1}^2  \frac{\gamma_{max}}{k^{1+\delta}}\|\nabla f(x_k)\| \|\nabla f_{\N_k}(x_k)\|.$$
Thus, we obtain 
\begin{eqnarray} \label{sl3a} 
 E[-\sum_{j=0}^2 \alpha_{k+j} \gamma_{k+j}\nabla f(x_k)^T \nabla f_{\N_k}(x_k) | \F_k] 
& = & - \frac{1}{k} \nabla f(x_k)^T E[\nabla f_{\N_k}(x_k) | \F_k]\nonumber \\
&& +  {   \frac{2 \gamma_{max}}{k^{1+\delta}}\|\nabla f(x_k)\| E[\|\nabla f_{\N_k}(x_k)\| | \F_k]}\label{bound0}\\
&\leq  &{    - \frac{1}{k} \|\nabla f(x_k)\|^2}\nonumber \\ 
&& +\frac{4 \gamma_{max} (\sqrt G +1) }{k^{1+\delta}} ({   \sqrt{G}}+\|\nabla f(x_k)\|^2), \label{bound1}
\end{eqnarray}
{   where the last inequality comes from \eqref{cor5}, considering both cases  $\|\nabla f(x_k)\|\leq 1, \|\nabla f(x_k)\|>1$.}
{ Finally, combining  \eqref{p1} with \eqref{p22}, \eqref{novo2} and \eqref{bound1} we obtain 
\begin{eqnarray*} \label{s4a} 
        E[f( x_{k+3}) -f_{low} | \F_k] & \leq &  f( x_{k}) -f_{low} -\frac{1}{k}\left( 1-\frac{ 4  \gamma_{max}(\sqrt G +1)+\tilde{C}_3 +{\color{magenta} L}\tilde{C}_2 /2}{k^{\delta}}\right)\|\nabla f(x_k)\|^2\\
        && + \frac{\tilde{C}_3 {   \sqrt{G}}+L\tilde{C}_2G /2 +{   4 (G+\sqrt{G})\gamma_{max}} }{k^2}
\end{eqnarray*}}
and conclude that there exists $\bar{k} \in K$ such that for  all $ k \geq \bar{k}, k \in K $ we have 
\begin{equation} \label{s5a} 
        E[f( x_{k+3}) -f_{low} | \F_k] \leq f( x_{k}) -f_{low} - \frac{1}{2 k} \|\nabla f(x_k)\|^2 +  \frac{\tilde{C}_3 {   \sqrt{G}}+{   L}\tilde{C}_2G {   /2}{   +4 (G+\sqrt{G})\gamma_{max}}}{k^2}.
\end{equation}
 Continuing as in {   Theorem \ref{Tuniss}} we obtain the result. }\end{prf} $\Box$

\subsection{Nonuniform sampling}
In this subsection we discuss the case where the sampling is not uniform.
 Let us denote the indices that form the subsample $ \N_k $
with $ i_1^k,\ldots, i^k_{{  S}},$ i.e.,
$  \N_k  = \left\{i_1^k,\ldots, i^k_{{  S}}\right\},  $ and let
 us define
$$ p_j^k:=P(k_\ell=j), \;  j=1,\ldots.N,$$
where $ k_1,\ldots,k_{{  S}}$ are i.i.d. random integers and $P(k_\ell=j)$ is the probability of the event $k_\ell=j$. We allow nonuniform sampling in the following sense. 

\begin{ass}[\bf{A4}] \label{A5}
There  exist constants $ M_3\geq 0$ and $\varepsilon>0$ such that
$$ |p_j^k - \frac{1}{N}| \leq \frac{M_3}{k^\varepsilon}, \; j=1,\ldots. N$$
\end{ass}

The above assumption actually means that the sampling we will employ converges to the uniform sampling although we allow it to be nonuniform  for arbitrary large number of initial iterations. For example, one can employ modified {   Adaptive Importance Sampling } algorithms \cite{liu}, for sample sizes of arbitrary size  as follows.
Starting with  set $ \pi_i^0, i=1,\ldots,N $, at a generic iteration $k\ge 1$ we define  the set of relevant gradient norm values,
$$ \pi^k_{i^{k}_j}=\|\nabla f_{i^{k-1}_j}(x_{k-1})\|, \; j = 1,\ldots,{   S}, $$
and {   $\pi^k_{i}=\pi^{k-1}_{i}$ for any $i \not \in  \N_{k-1} $.}
Then, we  calculate  the  probabilities as
\begin{equation} \label{eq14}
    p_j^k = \frac{1}{k^\varepsilon}\frac{\pi_j^k}{\sum_{i=1}^N \pi_i^k} +\left(1-\frac{1}{k^\varepsilon}\right)\frac{1}{N}, \; j=1,\ldots,N.
\end{equation}

 In order to compensate the bias that comes from nonuniform sampling, we state the  assumption on similarity of the local cost functions \cite{fs}. {   The assumption is satisfied for a number of important problems, namely the problems coming from machine learning like logistic regression and for quadratic cost function for bounded $x. $ }

\begin{ass}[\bf{A5}] \label{A4}
There exist constants $ M_1$ and $ M_2 $ such that for every $ x \in \mathbb{R}^n$ and every $ i \in \N $ there holds
$$\|\nabla f_i(x)\| \leq M_1 + M_2 \|\nabla f(x)\|.$$
\end{ass}


For sampling defined above we have the following result. Instead of \eqref{unb} we have sampling as in A\ref{A5} and instead of assumption A\ref{A2} we have A\ref{A4}. 
{   and we prove the a.s. convergence for
the modified version of the algorithm.}

\begin{te} \label{te:nonuniformModi}
    Let  A\ref{A1}, A\ref{A5} and A\ref{A4} hold.  Let $\{x_k\}$ be a sequence generated by  modified SLiSeS algorithm with steps S3'-S5' instead of S3-S5. Then, with probability 1, there holds
$\liminf_{k \to \infty}  \|\nabla f(x_k)\|=0.$
\end{te}
{   
\begin{prf}
Following the same steps as in Theorem \ref{Tuniss}, we obtain \eqref{eq8} {   for $ k $ such that $ mod(k,m)=0$ } and 
\begin{equation} \label{p1c1}
        E[f( x_{k+3}) -f_{low}| \F_k] \leq f( x_{k}) -f_{low}+ E[\nabla f(x_k)^T(x_{k+3}-x_k)| \F_k] + \frac{\tilde{C}_2 L}{2 k^2}  E[\|\nabla f_{\N_k}(x_k)\|^2 | \F_k].
\end{equation}
Then, to upper bound $E[-\sum_{j=0}^2 \alpha_{k+j} \gamma_{k+j}\nabla f(x_k)^T \nabla f_{\N_k}(x_k) | \F_k] $, notice that $\gamma_k=1/k$ and  $\alpha_k=1$, so they   are $\F_k$-measurable.
{   Then, we can proceed as in the proof of Theorem \ref{Tunimodi} and we obtain \eqref{bound0}.}
Given that the sampling is nonuniform, let us start estimating the gradient norm using the similarity assumption  A\ref{A5}. We have
\begin{equation} \label{eq23}
\|\nabla f_{\N_k}(x_k)\|  = {   \|\frac{1}{S} \sum_{j=1}^S \nabla f_{i_j^k}(x_k)\| }\leq M_1 + M_2 \|\nabla f(x_k)\|.
\end{equation}
 Moreover,
\begin{equation}
    \label{n11}
    E[\|\nabla f_{\N_k}(x_k)\|^2 | \F_k] \leq E[ 2 M^2_1 + 2 M^2_2 \|\nabla f(x_k)\|^2| \F_k]=2 M^2_1 + 2 M^2_2 \|\nabla f(x_k)\|^2.
\end{equation}
Furthermore, given that  $ E[\nabla f_{i}(x_k)|\F_k] = \sum_{j=1}^{N} p_j^k \nabla f_j(x_k), \; i=1,\ldots,N$ and
$$ E[\nabla f_{\N_k}(x_k)|\F_k]= \sum_{j=1}^{N} p_j^k \nabla f_j(x_k) - \nabla f(x_k) + \nabla f(x_k) =\sum_{j=1}^{N}(p_j^k - \frac{1}{N}) \nabla f_j(x_k) + \nabla f(x_k),  $$ using Assumptions A\ref{A4} and  A\ref{A5}, we get
\begin{eqnarray} \label{eq26}
  E[-\frac{1}{k} \nabla f(x_k)^T \nabla f_{\N_k}(x_k)|\F_k] & = & -\frac{1}{k} \sum_{j=1}^{N}(p_j^k - \frac{1}{N}) \nabla f(x_k)^T \nabla f_j(x_k) - \frac{1}{k} \|\nabla f(x_k)\|^2  \nonumber \\
&\leq &  - \frac{1}{k} \|\nabla f(x_k)\|^2 + \frac{1}{k}  \sum_{j=1}^{N}|p_j^k - \frac{1}{N}| \|\nabla f(x_k)\| \|\nabla f_j(x_k)\| \nonumber \\
& \leq &  - \frac{1}{k} \|\nabla f(x_k)\|^2  + \frac{1}{k} \sum_{j=1}^{N}\frac{M_3}{k^{\varepsilon}} \|\nabla f(x_k)\| \|\nabla f_j(x_k)\| \nonumber \\
& \leq &  - \frac{1}{k} \|\nabla f(x_k)\|^2  + \frac{M_3 }{k^{1+\varepsilon}} \sum_{j=1}^{N} \|\nabla f(x_k)\|(M_1 + M_2 \|\nabla f(x_k)\|) \nonumber \\
& \leq &  - \frac{1}{k} \|\nabla f(x_k)\|^2  + \frac{N M_1 M_3}{k^{1+\varepsilon}} \|\nabla f(x_k)\|+\\\nonumber
&& \frac{N M_2 M_3}{k^{1+\varepsilon}} \|\nabla f(x_k)\|^2.
\end{eqnarray}
Again, considering both cases:  $ \|\nabla f(x_k)\| \leq  1 $ and $ \|\nabla f(x_k)\| >  1 $ we conclude that
\begin{equation} \label{eq27}
 E[-\frac{1}{k} \nabla f(x_k)^T \nabla f_{\N_k}(x_k)|\F_k] \leq  - \frac{1}{k} \|\nabla f(x_k)\|^2 +\frac{\bar{C}_1}{k^{1+\varepsilon}} \|\nabla f(x_k)\|^2 + \frac{\bar{C}_1}{k^{1+\varepsilon}},
\end{equation}
with $ \bar{C}_1 =NM_3(M_1+M_2) $.
Now, using \eqref{p1c1},  \eqref{novo2},  \eqref{bound0} and \eqref{eq27} we obtain 
\begin{eqnarray} \label{eq28}
   E[f(x_{k+3})-f_{low} |\F_k] & \leq & f(x_k) - f_{low} -\frac{1}{k} \|\nabla f(x_k)\|^2+\frac{\bar{C}_3}{k^{1+\tau}}+\frac{\bar{C}_4}{k^{1+\tau}}\|\nabla f(x_k)\|^2,
\end{eqnarray}
where $\tau:=\min \{\delta, \varepsilon\}$,  $\bar{C}_3=2 \gamma_{max} M_1+ L \tilde C_2  M_1^2 + \bar C_1$ and $\bar{C}_4=\bar{C}_1+L \tilde C_2 M_2^2+2 \gamma_{max} M_2.$

Finally, we conclude that there exists $ \bar{k} $ such that for $ k \in K, k \geq \bar{k}$ we have
$$ E[f(x_{k+3})-f_{low} |\F_k]  \leq  f(x_k) - f_{low}  - \frac{1}{2k}\|\nabla f(x_k)\|^2 + \frac{\bar{C}_3}{k^{1+\tau}}, $$
and the statement follows analogously as in Theorem \ref{Tuniss}. $\Box$
\end{prf}}

 \subsection{{   P-L and Strongly convex functions}}

Notice that the previous result holds for nonconvex problems. However, if we assume that the objective function is strongly convex, we can prove that the sequence of outer iterations converges a.s. towards the unique solution of the original problem.

\begin{te} \label{TuniSC} Suppose that the assumptions of {   Theorem \ref{Tuniss}} hold {   and the function $f$ is strongly convex}.
Then, the sequence $\{x_{m k-1}\}_{k \in \mathbb{N}}$  converges to the solution of problem \eqref{prb} $x^*$ almost surely.
\end{te}

\begin{prf}
 Let us denote with $f^*$ the optimal value of problem \eqref{prb}.
Theorem \ref{Tuniss} implies that $\liminf_{k \to \infty} \|\nabla f(x_{m k-1})\|^2 =0$ a.s. The strong convexity further implies the existence of $K$ such that $\lim_{k \in K} x_{m k-1}=x^*$ a.s. Therefore, $\liminf_{k \to \infty}   f(x_{m k-1}) =f^*$ a.s. Moreover, using the same notation as in the proof of Theorem \ref{Tuniss} we conclude that Theorem~\ref{thm:rand_variable_conv} implies that $U_l$ converges a.s. which means that $\{f(x_{m k-1})\}_{k \in \mathbb{N}}$ converges a.s. Combining everything together, we conclude that the whole  sequence $\{f(x_{m k-1})\}_{k \in \mathbb{N}}$ converges to $f^*$ a.s. and due to strong convexity we conclude  $\lim_{k \to \infty} x_{m k-1}=x^*$ a.s.
$\Box$ \end{prf}


{   Relaxing the strong convexity assumption, we get the following result. 


\begin{te} \label{TuniPL} Suppose that the assumptions of  Theorem \ref{Tuniss} hold 
{   
and the  function $f$ satisfies the Polyak-Lojasiewicz (P-L) condition, i.e.,  for every $x \in \mathbb{R}^n$ there holds 
\begin{equation}\label{PL}
 \|\nabla f(x)\|^2 \geq 2 \mu (f(x)-f^*)
 \end{equation}
 where $\mu$ is a positive constant and $f^*$ is the optimal value of problem \eqref{prb}.
}
 Then, the sequence $\{f(x_{m k-1})\}_{k \in \mathbb{N}}$  converges to  $f^*$  almost surely.
\end{te}

\begin{prf}
Theorem \ref{Tuniss} implies that $\liminf_{k \to \infty} \|\nabla f(x_{m k-1})\|^2 =0$ a.s. and {    \eqref{PL} implies that $\liminf_{k \to \infty}  f(x_{m k-1}) =f^*$ a.s.. 
Proceeding as in the proof of Theorem \ref{TuniSC} it can be proved that $\lim_{k \to \infty}  f(x_{m k-1}) =f^*$ a.s..
}  
$\Box$ \end{prf}

{\bf{Remark.}} The two previous statements holds if instead of the assumptions in {\color{blue} } Theorem \ref{Tuniss} we impose the assumptions in   Theorems \ref{Tunimodi} or \ref{te:nonuniformModi}.

}

\section{Numerical results}

To give further insight into the proposed subsampling strategy, and to illustrate the practical performance of the SLiSeS algorithm,
we present the results of some numerical experiments.
  We consider two versions of our proposed algorithm: SLiSeS using uniform sampling (SLiSeS-UNI),  analyzed in Section 3.1,  and
 SLiSeS using {   Adaptive Importance Sampling}  (SLiSeS-AIS),  analyzed in Section 3.2 with probabilities defined through \eqref{eq14} with $\varepsilon=1$ and $\pi^0_i=1, i=1,...,N$. 
  {   Another variant, proposed and analyzed in Section 3.1.2, is to divide $c_k$ by  $k^{1+\delta}$, where $\delta>0$  is a small  given number.
   We will denote this variant as the Modified SLiSeS method.  All these versions will be compared with several recently
   proposed schemes in the literature for solving (\ref{prb}). }

Although the SLiSeS algorithm can be applied more generally, in our experiments we
focus on strongly convex quadratic functions and also on L2-regularized logistic regression problems,
as problems that satisfy our {   assumption A\ref{A4}.} The considered quadratic functions, for $1\leq i\leq N$, are given by
$$ f_i(x) = \frac{1}{2}(x-b_i)^{\top} A_i (x-b_i), $$
 where  $b_i \in \rn$ and  the symmetric positive definite matrix $A_i \in \r^{n\times n}$ are obtained as in \cite{jkk},
  i.e.,  vectors $b_i$  are extracted from the Uniform distribution on [1,31], independently from each other. Matrices
   $A_i$ are of the form $A_i = Q_i D_i Q_i^{\top}$, where $D_i$ is a diagonal matrix with Uniform distribution on [1,101] and
   $Q_i$ is the matrix of orthonormal eigenvectors of $\frac{1}{2}(C_i + C_i^{\top})$, and the matrix $C_i$ has components
   drawn independently from the standard Normal distribution. For the logistic regression problems, for $1\leq i\leq N$, the function $f_i$ is given by
   $$ f_i(x) = \log(1 + \exp(-b_i(a_i^{\top} x))) +  \frac{\lambda}{2}\|x\|_2^{2}, $$
   where the vectors $a_i\in \rn$ and the labels $b_i \in \{-1,1\}$ are given, and
    the regularization parameter $\lambda$ is set to  $10^{-4}$. The vectors $a_i$ and the labels $b_i$ are
    obtained from three real data sets: Voice (LSVT Voice Rehabilitation) \cite{Voice} ($n=309$, $N=126$),
  Park (Parkinson's Disease Classification) \cite{Park} ($n=754$, $N=756$), and the  Cina0
    econometrics dataset ($n=132$, $N=16033$)  which was downloaded from the CINA website\footnote{http://www.causality.inf.ethz.ch/data/CINA.html}.

  In all our experiments,  we report the computational cost measured by the number of  function evaluations versus
 the observed decrease in the objective function. To calculate the cumulative number of function
 evaluations   
 we add $S=|\N_k|$ each time we evaluate the function $f_{\N_k}$, 
 	and the gradient evaluation comes for free; see, e.g., \cite[p. 27]{bkm}.
 All runs are terminated when a maximum number of iterations ($maxiter$) is reached.  The required input parameters in the SLiSeS algorithms
 are fixed as follows: $x_0=0$, $\eta=10^{-4}$, $t_k=1/2^k$,  $\gamma_{min} = 10^{-8}$,  {   $\gamma_{max}=10^{8}$ and $\tilde \gamma=1/\|g_0\|$ for the Modified SLiSeS.}
 } Concerning
 the sample size, unless otherwise specified, we set $S=|\N_k|=1$ for all $k$.


\begin{figure}[htbp]
 \includegraphics[height=6.5cm, width=17cm]{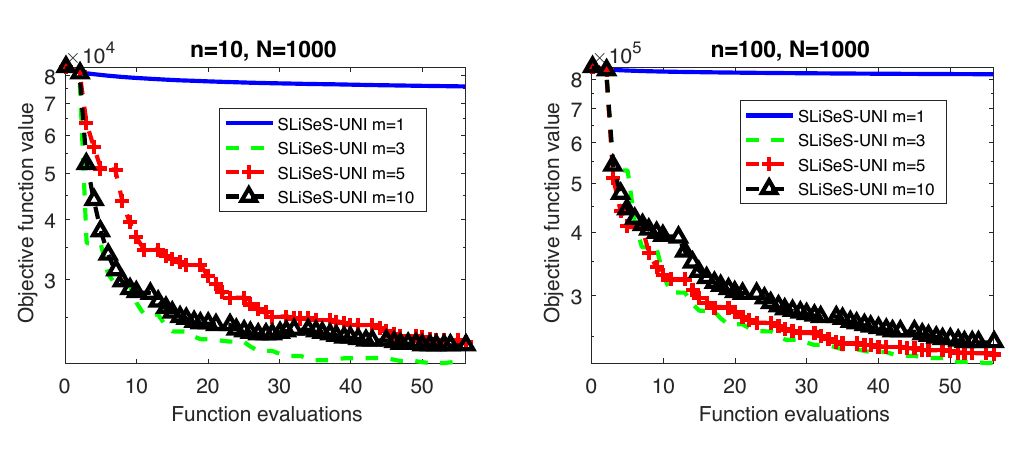}
	\caption{{\footnotesize{Performance of SLiSeS-UNI on strictly convex quadratics for  $m=1,3,5,10$ (inner iterations), $N=1000$, $maxiter = 50$,
				and $n=10$ (left) or $n=100$ (right). }} }
	\label{m1-10}
\end{figure}

\begin{figure}[htbp]
	\includegraphics[height=6.5cm, width=17cm]{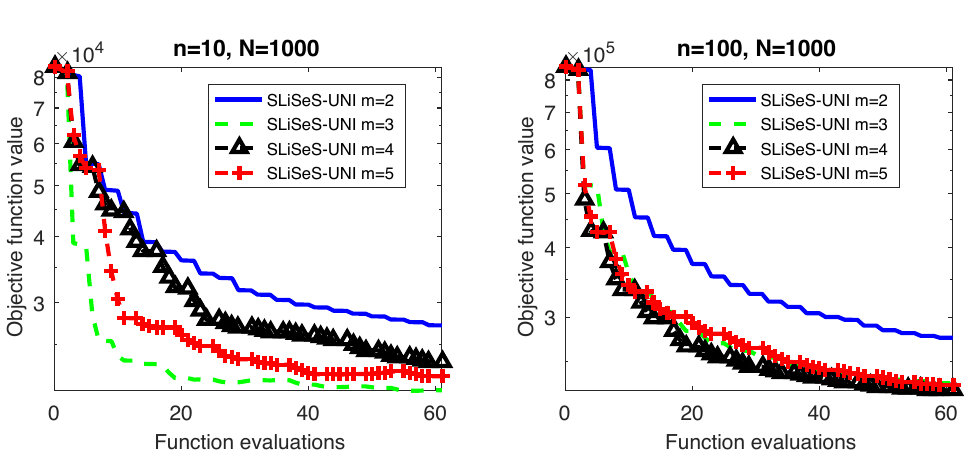}
	\caption{{\footnotesize{Performance of SLiSeS-UNI on strictly convex quadratics for  $m=2,3,4,5$ (inner iterations),  $N=1000$, $maxiter = 50$,
				and $n=10$ (left) or $n=100$ (right). }} }
	\label{m2-5}
\end{figure}

\begin{figure}[htbp]
 	\hspace*{-0.1in}
	\includegraphics[height=10.0cm, width=9cm]{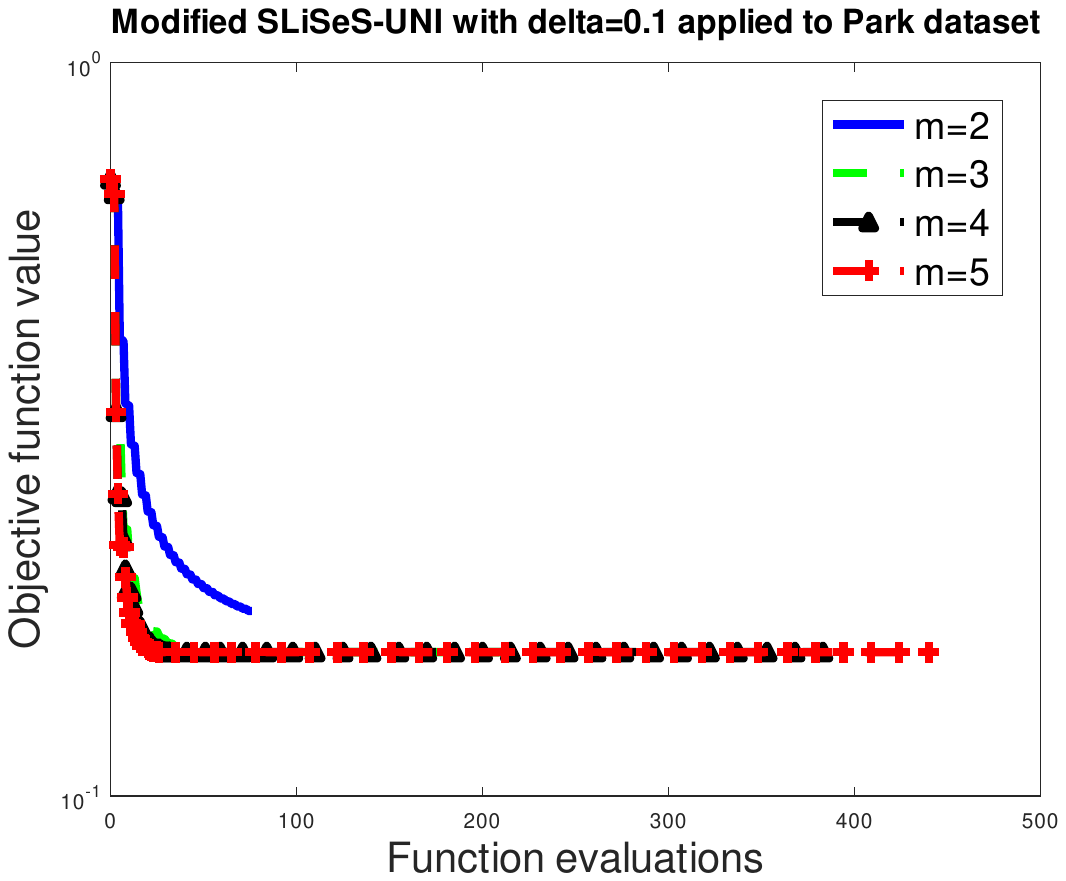}
   	\hspace*{-0.3in}
    \includegraphics[height=10.0cm, width=9cm]{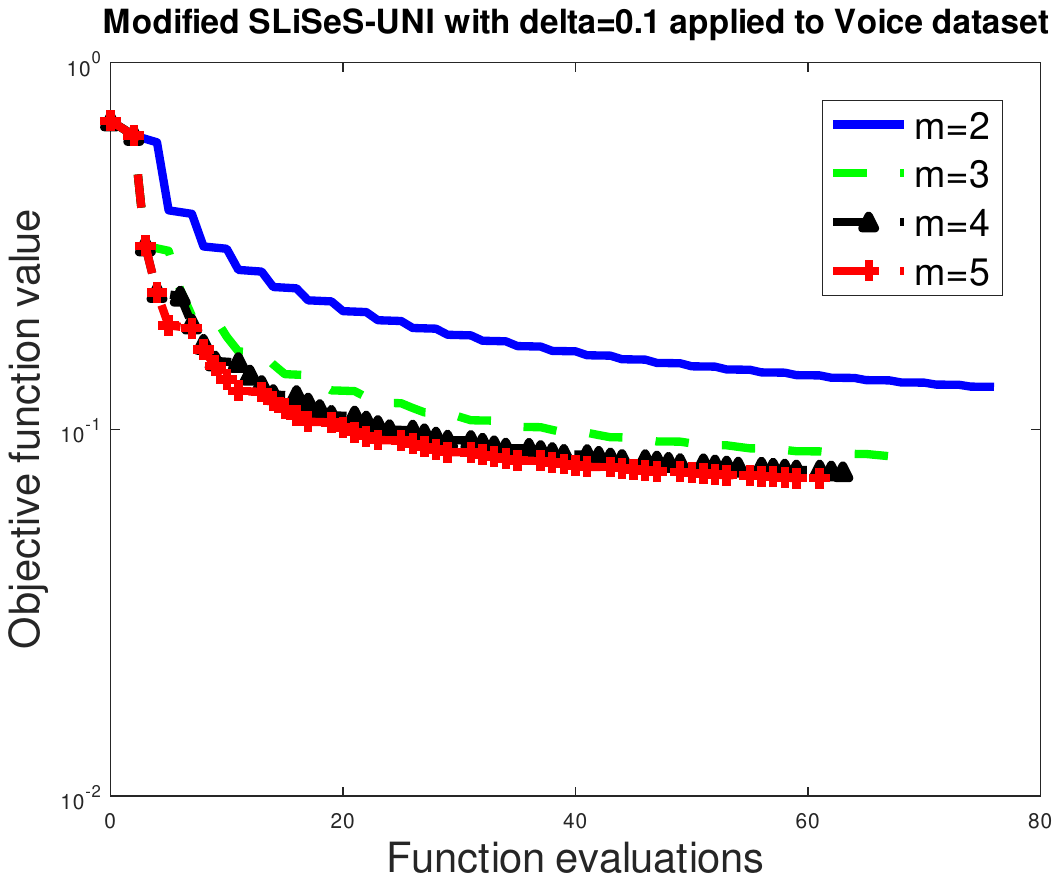}
	\vspace*{-0.8in}
	\caption{{\footnotesize{Performance of {   Modified SLiSeS-UNI with $ \delta=0.1 $}  on the Park dataset (left) and on the Voice dataset (right)
				for  $m=2,3,4,5$ (inner iterations), and $maxiter = 50$.  }} }
	\label{VPm2-5}
\end{figure}

We  begin by exploring the behavior of  Algorithm SLiSeS-UNI {   and Modified SLiSeS-UNI with $\delta=0.1$} for different values of the number of internal iterations $m$,
  in which we keep the same sample. 
  Let us start by considering  strictly convex quadratic functions. In Figure  \ref{m1-10} we report
   the performance of Algorithm SLiSeS for several values of $m$ from 1 to 10, $N=1000$, $maxiter = 50$, and two different values of $n$: $10$ and
    $100$. We can observe that, for both values of $n$, the choice $m=1$ (i.e., the standard strategy of changing the sample at each iteration)
     produces a much worse behavior than any of the other 3 choices of $m>1$. We also note, for this particular experiment as well
       as many others we have run,  that the best behavior is observed when $2\leq m \leq 5$. In Figures \ref{m2-5} and \ref{VPm2-5}
        we focus our attention on that interval, for the strictly convex quadratic case and also  for some logistic
       regression problems (Voice and Park datasets).
       {   For these latter test problems we report results with the modified version as we observe a general better behaviour of the modified versions with respect to the unmodified counterpart on this class of problems, as it will be shown in Figure 6.}
   We note that for $m=3$, $m=4$, $m=5$ the performance is better than for $m=2$, and also that none of the choices 3, 4, or 5 shows a clear advantage
          over the other two. However, in the average, the value $m=3$ seems to be the one that shows a better performance. 

\begin{figure}[htbp]
  	 \hspace*{-0.1in}
   \includegraphics[height=20cm, width=16cm]{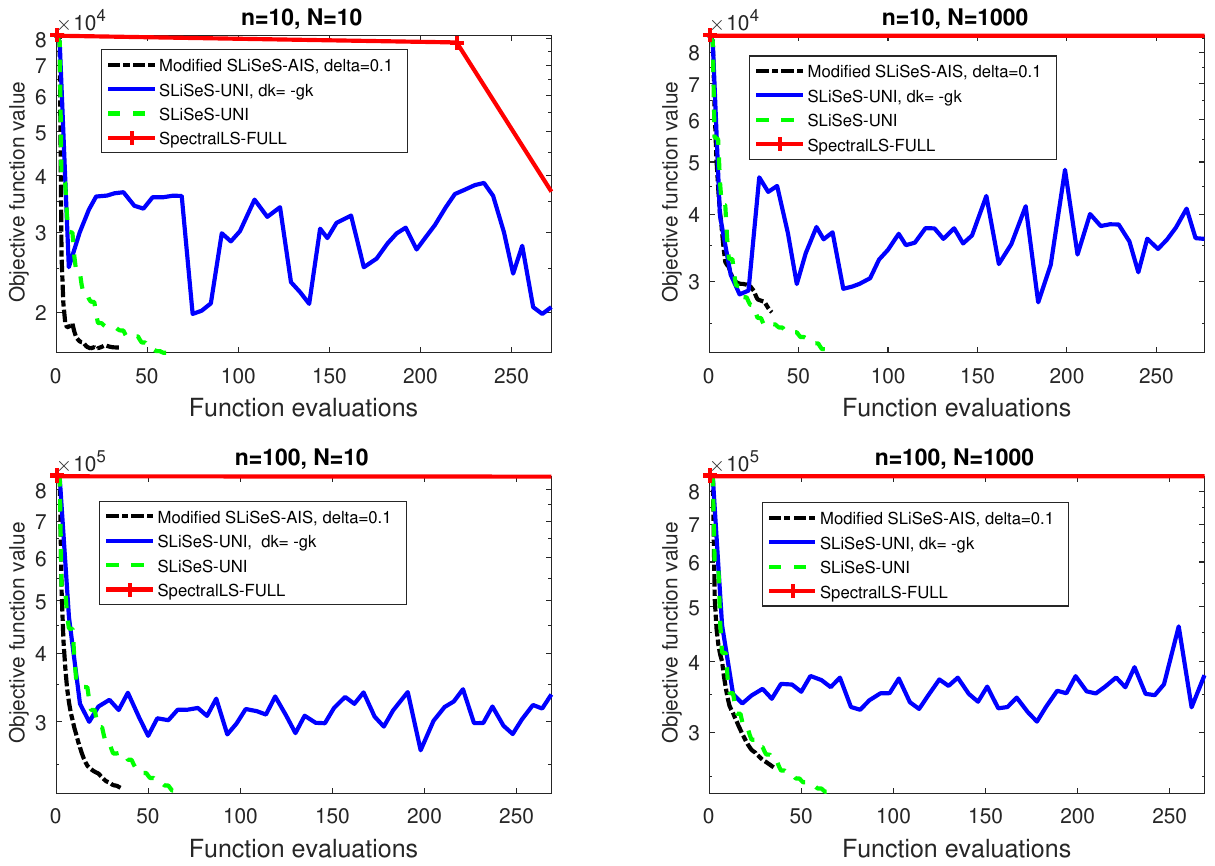}
 	\vspace*{-1.6in}
 	\caption{{\footnotesize{Performance of  SLiSeS-UNI ($m=3$), SpectralLS-Full, {  Modified SLiSeS-AIS
 		(for $\delta=0.1$)} and {   SLiSeS-UNI with $d_k=-g_k$,} for
 				strictly convex quadratic functions, different combinations of $n$ and $N$, and $maxiter = 50$. }} }
 	\label{UniNonFullq}
 \end{figure}

For our second set of experiments, we compare   SLiSeS-UNI {   and Modified SLiSeS-AIS with $ \delta= 0.1$
} with the full sampled
 version of the SLiSeS algorithm, i.e., when  $S=|\N_k|=N$ for all $k$. Notice that this so-called SpectralLS-Full
 method matches  the spectral gradient method globalized with the line search procedure described in Algorithm 2,
  for which global convergence is guaranteed; see \cite{bmr14} and references therein. Hence, in the SpectralLS-Full method,
  there is no need to divide by $k$ at Step S4 of Algorithm 1. {   To illustrate the benefits of using spectral-BB
 steplegths in our approach, we also consider the same Algorithm 1 with two additional parameter settings. In the first one we fix $ \gamma_k=1$, so $ d_k = - g_k,$ while in the second one we fix $ \gamma_{min} = \gamma_{max} =1, $ and so $ d_k = - (1/k) g_k. $ {   In both of these latter procedures we  inhibit the employment of BB steplenghts in the SLiSeS framework, and}   
neither of these two settings is competitive  with SLiSeS employing  the spectral-BB step sizes. } In Figure  \ref{UniNonFullq}  we report
  the performance of Algorithms  SLiSeS-UNI ($m=3$),  SpectralLS-Full, {   SLiSeS-UNI with  $ d_k =-g_k $ and Modified SLiSeS-AIS with $ \delta= 0.1$  when applied to
  strictly convex quadratics for $maxiter = 50$ and  different combinations of $n$ and $N$. The case SLiSeS-UNI with $\gamma_{min} = \gamma_{max} =1 $ is not shown in the graphs as it is quite the same as $ \gamma_k =1; $ i.e SLiSeS-UNI with  $ d_k =-g_k $.}
   We also note that   for any combination of $n$ and $N$, the SpectralLS-Full method requires a much higher number of function evaluations, so it is not competitive.

\begin{figure}[htbp]
 \includegraphics[height=7cm, width=12cm]{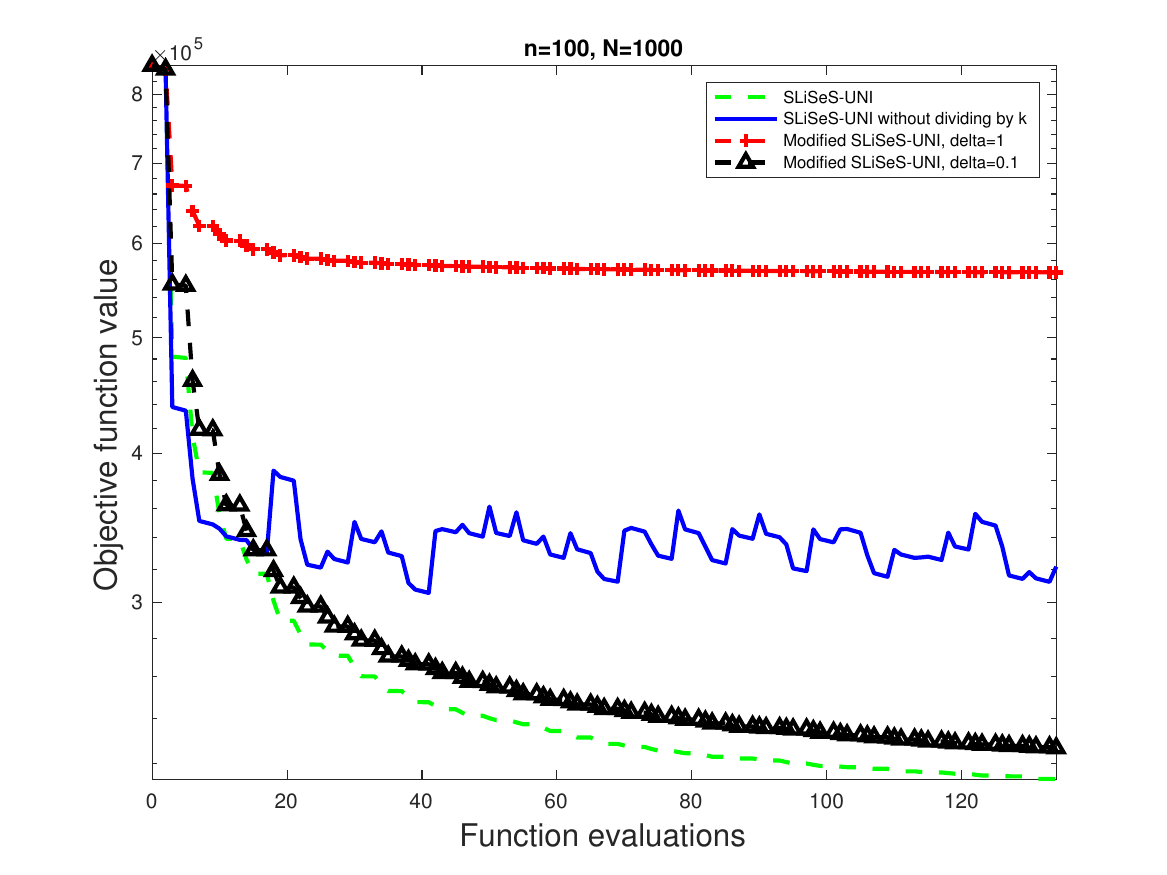}
\caption{{\footnotesize{Performance of the  SLiSeS-UNI method without dividing by $k$ at Step S4, and the Modified SLiSeS-UNI method (for
				$\delta=0.1$ and $\delta=1$) when  compared with the original versions of   SLiSeS-UNI  for $m=3$, $maxiter=100$,
				$n=100$,  $N=1000$ and strictly convex quadratic functions.}} }
	\label{ModifiedV}
\end{figure}

 Our convergence analysis in Section 3 was developed using the damping strategy of dividing by $k$ the spectral step length
 $c_k$  at Step S4 of the SLiSeS algorithm. Nevertheless, in order to assess the benefits of that damping strategy,
  in our next set of experiments we explore the practical performance of the  SLiSeS algorithm when $c_k$ is not divided by $k$ at Step S4.
 In Figure  \ref{ModifiedV}  we show
  the performance of 
  SLiSeS-UNI ($m=3$) with the two mentioned variants,
    when applied to strictly convex quadratics for $maxiter = 100$, $n=100$, $N=1000$, and two different values of $\delta$ ($0.1$ and 1).
 
\begin{figure}[htbp]
	\hspace*{-0.1in}\includegraphics[height=10.0cm, width=9.0cm]{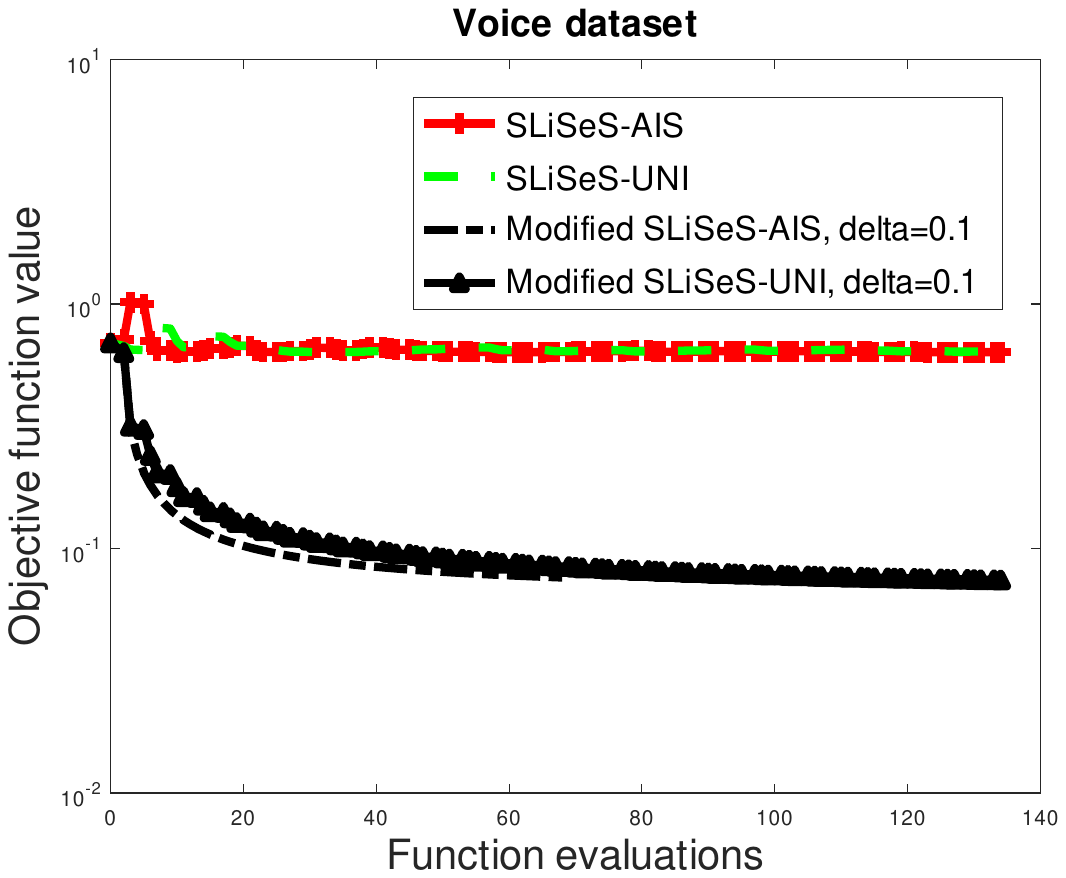}
	\hspace*{-0.3in}\includegraphics[height=10.0cm, width=9.0cm]{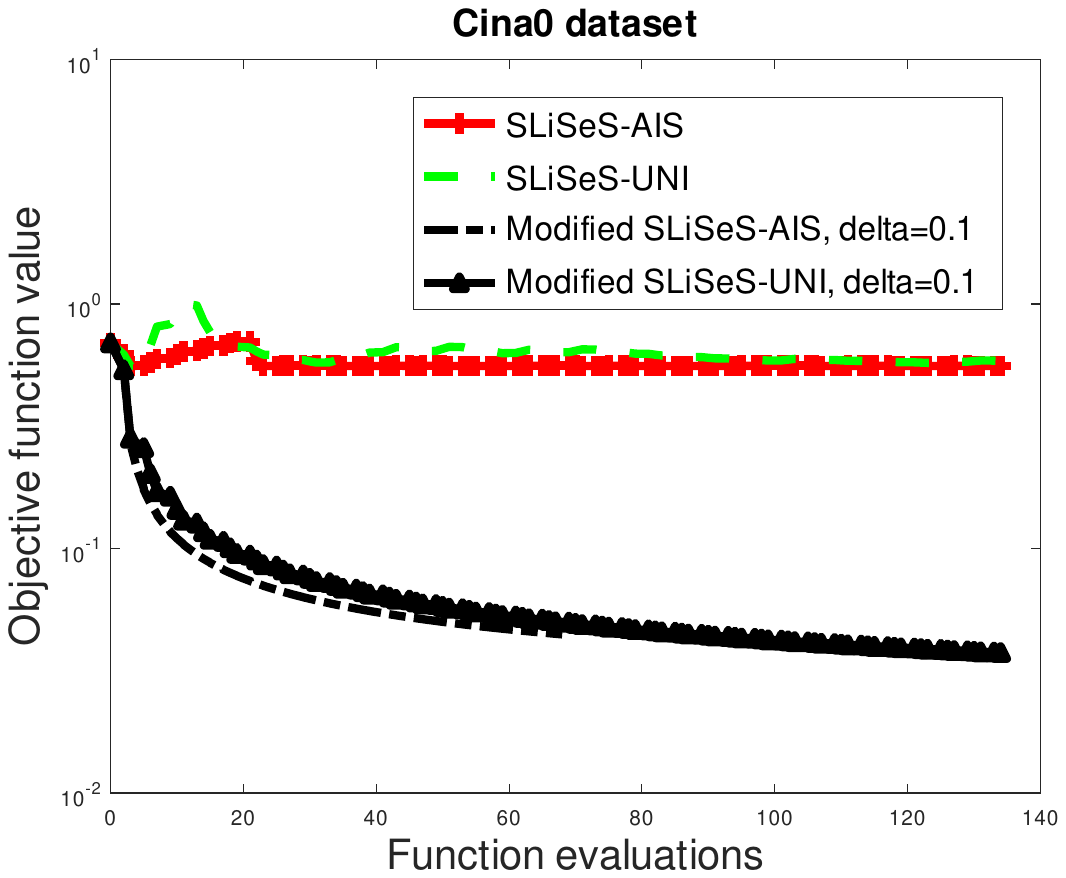}
	 \vspace*{-65pt}
 \caption{{\footnotesize{{   Performance of the SLiSeS methods (AIS and UNI) and their modified versions  (with $\delta=0.1$) for $m=3$, on the  Voice  dataset (left),  and on the Cina0 dataset (right).}} }}
	\label{newfigure}
\end{figure}

     It can  be noted that in the case of strictly
      convex quadratics the version that does not divide by $k$ produces unstable behavior. Moreover, for many of the experiments we have tried,
       in the case of logistic regression functions that variant does not converge. {   We can also observe 
       	that the damping strategy using $\delta=0.1$ shows a better performance than using $\delta=1$. In
       Figure \ref{newfigure} we present the advantages of using the damping strategy with $\delta=0.1$ in the
        case of two logistic regression functions. These last two figures illustrate that the damping strategy
        not only guarantees convergence but also adds practical advantages to the SLiSeS algorithm.}

\begin{figure}[htbp]
	\hspace*{-0.25in}\includegraphics[height=10.0cm, width=9.0cm]{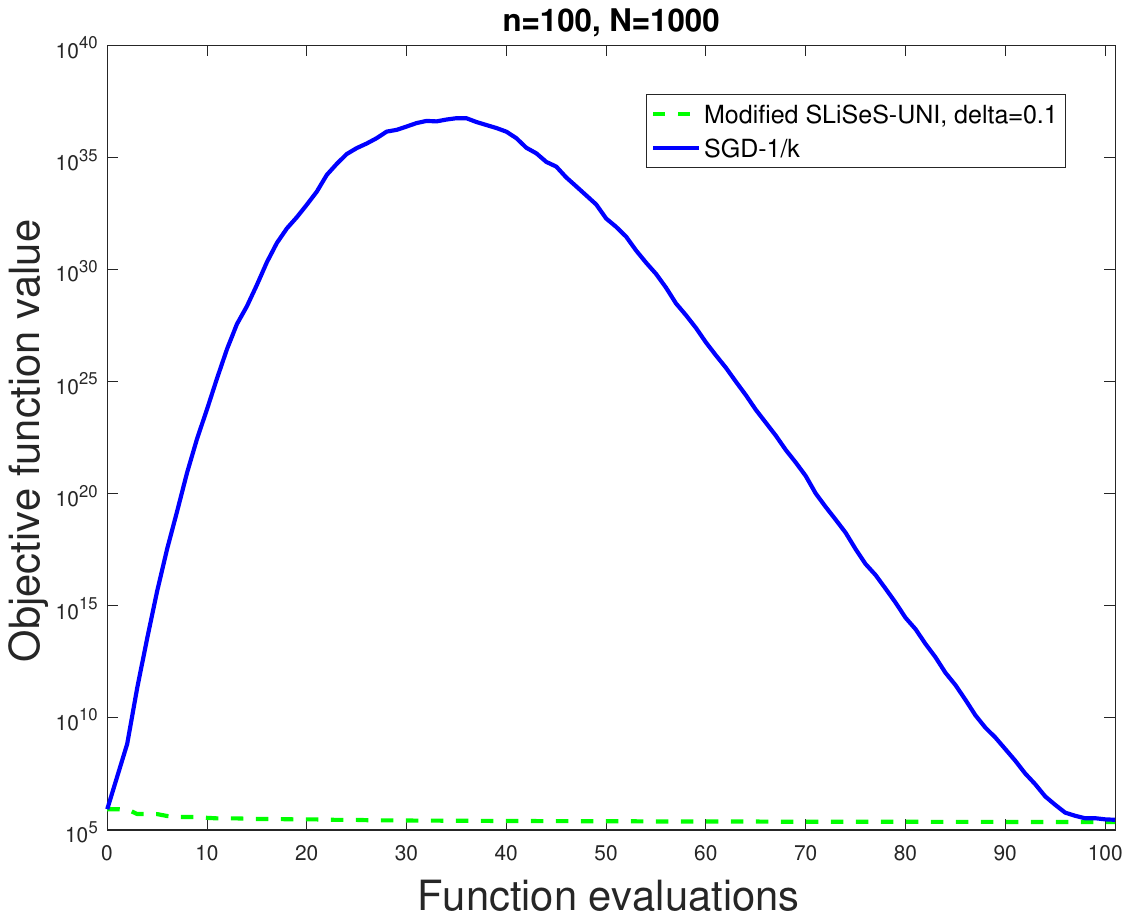}
 \hspace*{-0.3in}	\includegraphics[height=10.0cm, width=9.0cm]{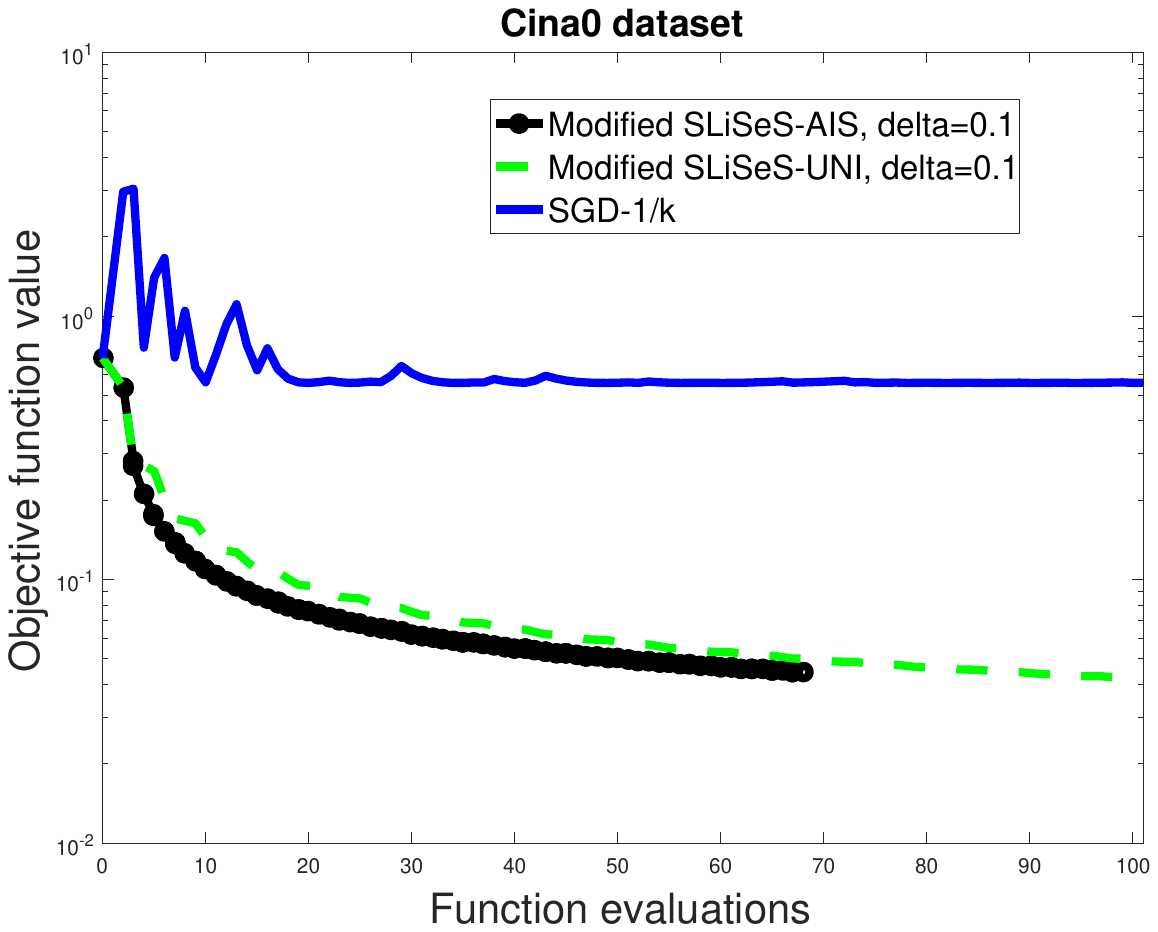}
	   \vspace*{-50pt}
\caption{{\footnotesize{{   Performance of the modified SLiSeS-UNI  method ($m=3$ and $\delta=0.1$) and the SGD method with step length $1/k$ ($maxiter=100$) on a strictly convex quadratic when $n=100$ and 
			$N=1000$  (left),  and performance of the modified versions of SLiSeS ($\delta=0.1$) and the SGD method with step length $1/k$ ($maxiter=100$) on the Cina0 dataset (right).}}} }
	\label{SLiSeSvsSGD}
\end{figure}

\begin{figure}[htbp]
 \hspace*{-0.25in}\includegraphics[height=10cm, width=9.5cm]{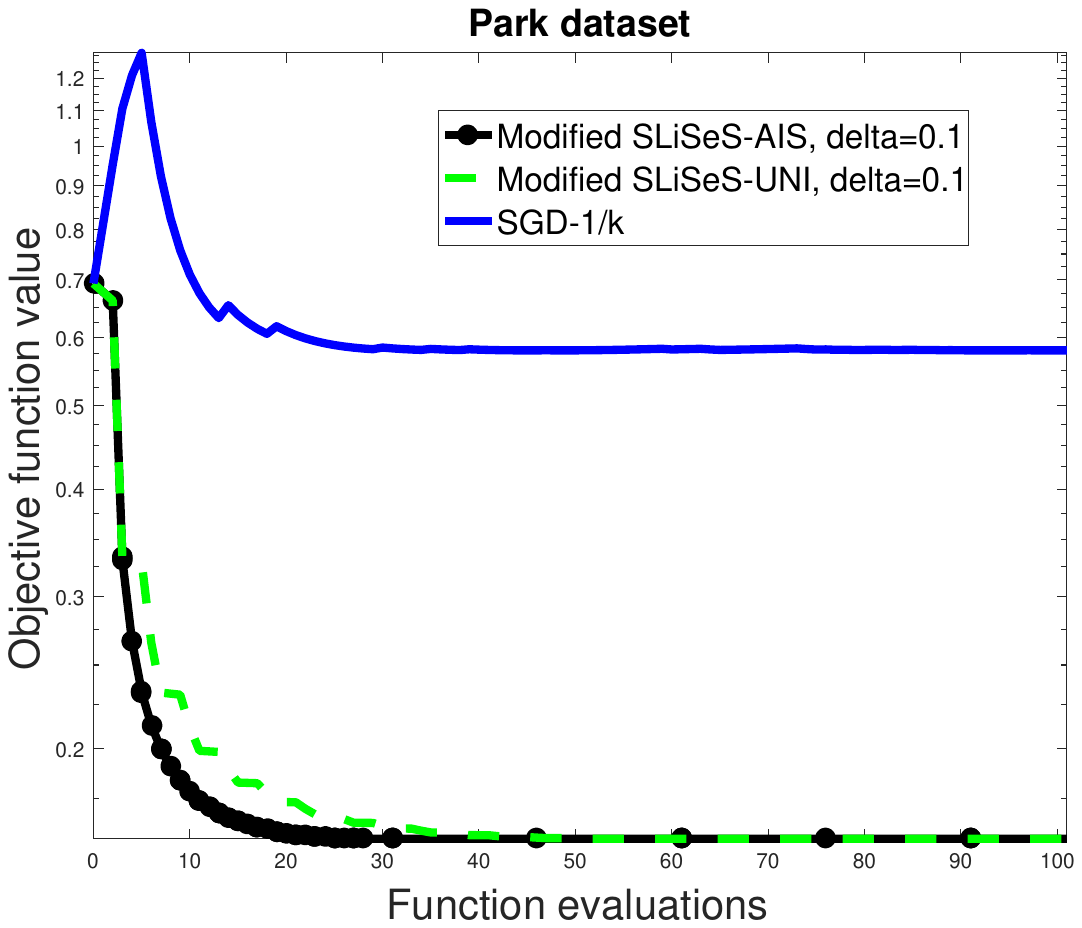}
 \hspace*{-0.2in}\includegraphics[height=10cm, width=8.4cm]{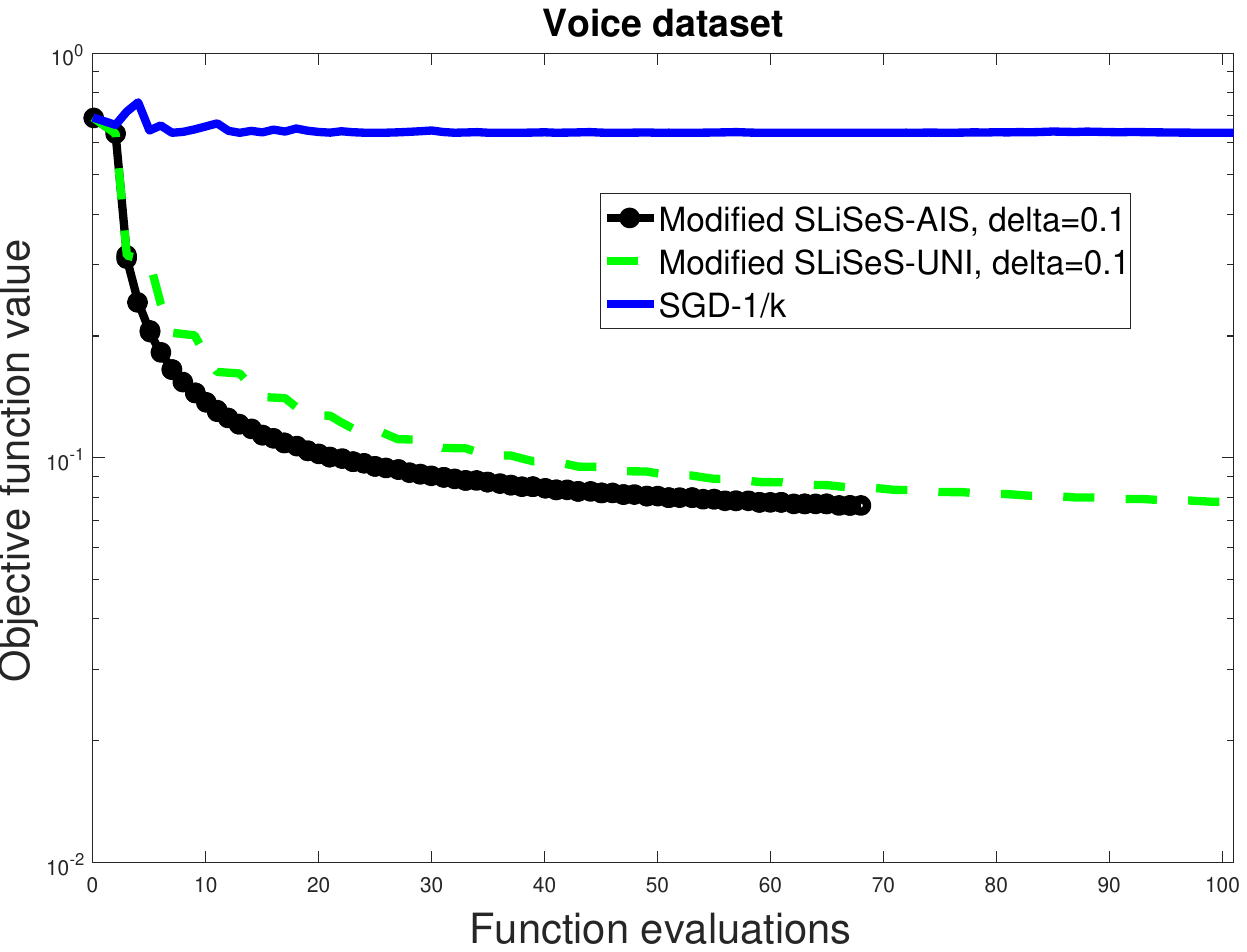}
 \vspace*{-0.8in}
	\caption{{\footnotesize{{   Performance of the modifed SLiSeS methods ($m=3$ and $\delta=0.1$) and the SGD method with step length $1/k$, for $maxiter=100$, on the Park dataset (left) and the Voice dataset (right).}}} }
	\label{SLiSeSvsSGD2}
\end{figure}

Finally, in the last set of experiments we are interested in comparing the  SLiSeS method with the natural competitors that exist in the literature
 to solve problem (\ref{prb}). The so-called Stochastic Gradient Descent (SGD) continues to receive significant attention, and is still
  a widely used technique for machine learning applications; see, e.g., \cite{Bottou_et_al} and references therein.
 The SGD method can be interpreted as a simplified case of the SLiSeS algorithm where the step length $\gamma_k$
 is obtained either as a sufficiently small constant  to guarantee convergence, or in a diminishing way
 such that $\sum \gamma_k = \infty$ and $\sum \gamma_k^2 < \infty$ to guarantee convergence
  (the most common choice is $\gamma_k = 1/k$). In either case, the SGD method does not employ a line search strategy.
  In Figure  \ref{SLiSeSvsSGD}  we report the performance of {   the modified variants of SLiSeS 
  ($m=3$ and $\delta=0.1$)} and the SGD method with step length $1/k$, for $maxiter=100$, 
  when applied to a  a strictly convex quadratic when $n=100$ and $N=1000$, and also to the Cina0 dataset. 
  In Figure  \ref{SLiSeSvsSGD2} we show  the performance of the same {   three} methods when applied to the Park dataset and also to the Voice dataset.
  The results observed in both figures are indicative of the known weaknesses of the SGD method, and also illustrate that the {   modified SLiSeS Algorithms}  exhibits better convergence behavior.

 More recently, the Barzilai and Borwein spectral step length has been used in connection with stochastic gradient ideas in the so-called  SVRG-BB
 family of methods; see  \cite{tmdq}. The SVRG-BB uses a full gradient (i.e., $S=N$)  every $p$ iterations (in \cite{tmdq}, it is recommended to use
   $p=2n$), and it depends on a certain parameter $\eta_0>0$ (in \cite{tmdq}, it is recommended to set $\eta_0=0.01$).  Two more members of
    that family are proposed and discussed in \cite{tmdq}, the so-called SGD-BB and the SGD-BB with a
    specialized smoothing technique (both of them  avoid calculating full gradients).
    These last two options depend on four parameters: $\beta>0$, $\eta_0>0$, $\eta_1>0$, and $p$.
    (in \cite{tmdq}, it is recommended to set $\eta_0=\eta_1= 0.01$, $p=n$, and $\beta=1/p$). In \cite{tmdq}, convergence properties of the
    SVRG-BB family are established for strongly convex functions.
    
    \begin{figure}[htbp]
    	\hspace*{-0.15in}\includegraphics[height=10cm, width=9.0cm]{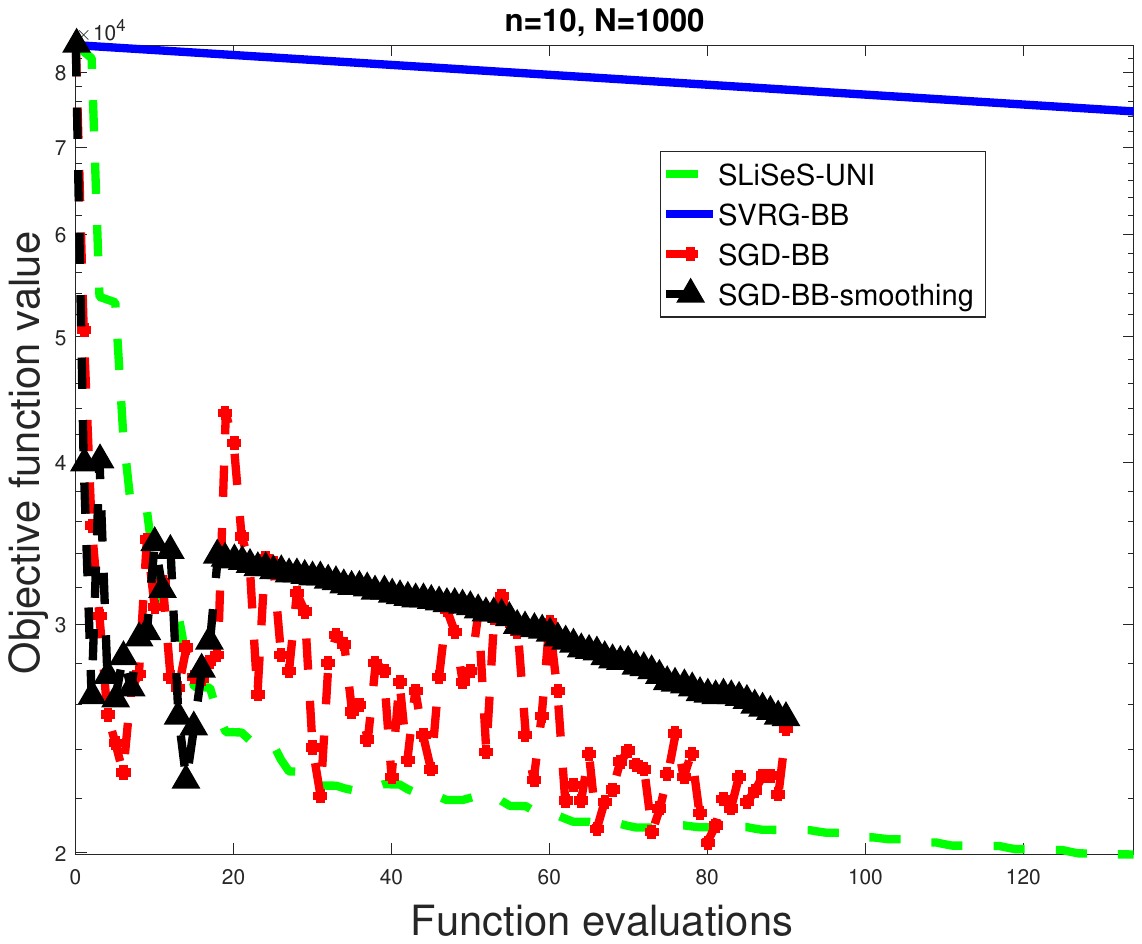}
    	\hspace*{-0.15in}\includegraphics[height=10cm, width=9.0cm]{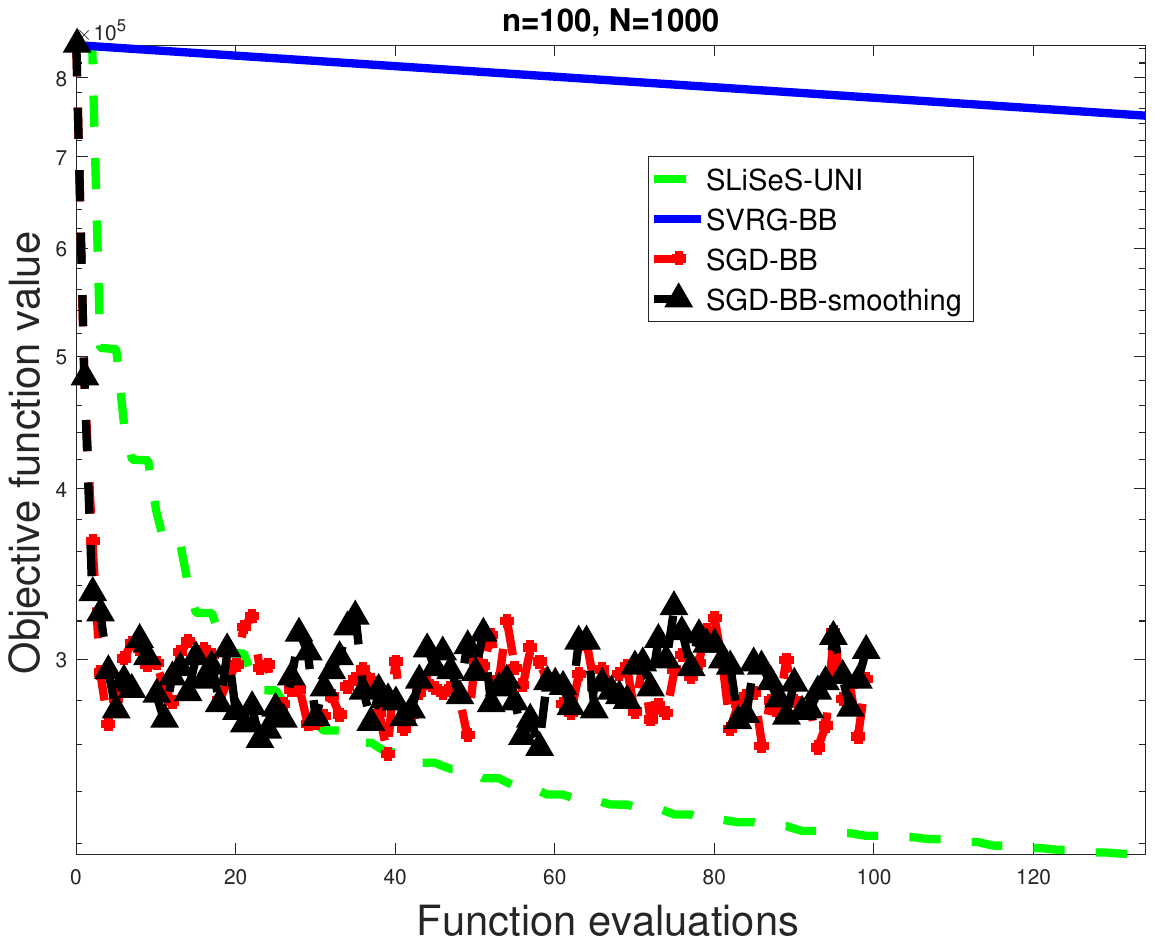}
    	\vspace*{-0.9in}
    	\caption{{\footnotesize{Performance  of the  SLiSeS-UNI ($m=3$) on a strictly convex quadratic
    				when compared with SVRG-BB ($p=2n$ and $\eta_0>0.01$),  SGD-BB ($p=n$, $\beta=1/p$, and $\eta_0=\eta_1=0.01$), and  SGD-BB with a  smoothing technique
    				(same parameters as SGD-BB) for  $maxiter=100$, $N=1000$, and $n=10$ (left) and $n=100$ (right).}} }
    	\label{svrgbb}
    \end{figure}

     In Figure  \ref{svrgbb} we report, for strictly convex quadratics, the performance of the SLiSeS-UNI ($m=3$) method and the three described variants of the SVRG-BB family of methods:  SVRG-BB,  SGD-BB, and  SGD-BB with a  smoothing technique, when $maxiter=100$, $N=1000$, and $n=10$ (left) and $n=100$ (right). We note that the SVRG-BB method has a much worse performance than all the other
     competitors, and this is clearly due to the use of a full gradient every $p$ iterations. We note that the best performance is obtained
      by the SLiSeS-UNI ($m=3$) method. We also note that for small values of $n$ the
      smoothing technique adds some benefit to the SGD-BB method, but once $n$ increases the smoothing effect tends to vanish.

\section{Conclusions}

Stochastic algorithms based on the negative gradient direction remain one of the most convenient options for solving large-scale finite-sum optimization problems. Among these, the well-known SGD method and its variants present attractive theoretical advantages since they allow establishing worst-case complexity results. However, from a practical point of view, the SGD method has some known weaknesses and can exhibit very slow convergence.

Recently, focusing on practical behavior, some  ideas have incorporated the use of non-monotonic step lengths such as the Barzilai and Borwein (BB)
 spectral one. These ideas have recently been proposed and analyzed in the literature, either avoiding full gradient  evaluations altogether or
 by allowing full gradient evaluations at some periodic iterations. For all these gradient-type methods, including the SGD method and its variants, a common denominator has been to keep the tradition of changing the sample at each iteration.

In this work, following the line of incorporating the BB spectral step length, and accepting its non-monotonic behavior, we proposed and analyzed a novel
strategy: keeping the sample unchanged for {   a prefixed number of } consecutive iterations before sampling again.
The motivation behind this strategy is to promote the widely observed sweeping-spectrum behavior associated with spectral BB step lengths.
 The variants that emerge from this strategy have been
analyzed, and  convergence results have been obtained for non-convex smooth functions, for both uniform and nonuniform sampling. To illustrate and deepen
 this sampling strategy, we show a variety of numerical experiments that include comparisons with natural candidates from the literature. Our most
 important conclusion is that this strategy, based on the experiments shown, is promising and effectively improves the behavior of gradient-type methods.
   It deserves to be further analyzed with specific applications in mind. \\ [2mm]

\noindent
 {\bf Acknowledgements.}  We would like to thank an anonymous reviewer for  constructive comments 
 and suggestions that helped us to improve the final version of this paper.
 \\ [2mm]

\noindent
 {\bf Conflicts of interest.} No potential conflict of interest was reported by the authors. \\ [2mm]

\noindent
{\bf Funding.}
 The first author acknowledges financial support received  by the INdAM GNCS  and 
by PNRR - Missione 4 Istruzione e Ricerca - Componente C2 Investimento 1.1, Fondo per il Programma Nazionale di Ricerca e Progetti di Rilevante Interesse Nazionale (PRIN) funded by the European Commission under the NextGeneration EU programme, project ``Advanced optimization METhods for automated central veIn Sign detection in multiple sclerosis from magneTic resonAnce imaging (AMETISTA)'',  code: P2022J9SNP, MUR D.D. financing decree n. 1379 of 1st September 2023 (CUP E53D23017980001), project ``Numerical Optimization with Adaptive Accuracy and Applications to Machine Learning'',  code: 2022N3ZNAX
 MUR D.D. financing decree n. 973 of 30th June 2023 (CUP B53D23012670006).
 The second and  third authors were financially supported by the  Science Fund of the Republic of Serbia, Grant no. 7359, Project LASCADO.
 The fourth author was financially supported by Funda\c c\~ao para a Ci\^encia e a Tecnologia (Portuguese Foundation for Science and Technology)
		under the scope of the projects UIDB/MAT/00297/2020 (doi.org/10.54499/UIDB/00297/2020), and UIDP/MAT/00297/2020   (doi.org/10.54499/UIDP/00297/2020) 
		(Centro de Matem\'atica e Aplica\c{c}\~oes).  \\  [2mm]

 \noindent {\bf Data availability.}  The codes and datasets generated during and/or analyzed during the current study are
 available from the authors on reasonable request. \\ [2mm]

\end{document}